\documentclass[11pt, twoside]{article}

\usepackage{amssymb}
\usepackage{amsmath}

\allowdisplaybreaks

\def\rr{{\mathbb R}}
\def\rn{{{\rr}^n}}
\def\rrm{{{\rr}^m}}
\def\rnm{{{\rr}^n\times\rr^m}}
\def\rnt{{{\rr}^{n_1}\times\rr^{n_2}}}
\def\zz{{\mathbb Z}}
\def\nn{{\mathbb N}}
\def\cc{{\mathbb C}}
\def\cm{{\cal M}}
\def\cs{{\cal S}}

\def\ca{{\mathcal A}}
\def\cb{{\mathcal B}}
\def\cbq{{\cb_q}}
\def\cd{{\mathcal D}}
\def\crz{{\mathcal R}}
\def\boz{\Omega}
\def\fz{\infty}
\def\az{\alpha}
\def\supp{{\rm{\,supp\,}}}

\def\loc{{\rm{\,loc\,}}}
\def\lip{{\rm{\,Lip\,}}}
\def\lz{\lambda}
\def\dz{\delta}

\def\ez{\epsilon}
\def\bz{\beta}
\def\vz{\varphi}
\def\gz{{\gamma}}

\def\sz{\sigma}
\def\boz{\Omega}

\def\ls{{\lesssim}}
\def\laz{{\lfloor}}
\def\raz{{\rfloor}}

\def\wz{\widetilde}
\def\hs{\hspace{0.5cm}}

\def\laa{{\lip(\az_1,\,\az_2;\,\rnm)}}

\def\lp{{L^p(\rnm)}}
\def\lq{{L^q(\rnm)}}

\def\hp{{H^p(\rnm)}}

\def\dsum{\displaystyle\sum}
\def\diam{{\rm{\,diam\,}}}
\def\dint{\displaystyle\int}

\def\dfrac{\displaystyle\frac}
\def\dsup{\displaystyle\sup}

\def\r{\right}
\def\lf{\left}

\newtheorem{thm}{\sc Theorem}[section]
\newtheorem{lem}{\sc Lemma}[section]

\newtheorem{rem}{\sc Remark}[section]
\newtheorem{cor}{\sc Corollary}[section]
\newtheorem{defn}{\sc Definition}[section]
\newtheorem{pf}{\it Proof.}

\begin{document}

\baselineskip=15pt
\renewcommand{\arraystretch}{2}
\arraycolsep=1pt

\title{{\vspace{-4cm}\small\hfill\bf J. Math. Soc. Japan, to appear}\\
\vspace{4cm}\bf\Large Boundedness of Sublinear Operators
on Product Hardy Spaces and Its
Application\footnotetext{\hspace{-0.2cm}2000 {\it Mathematics Subject
Classification:} Primary: 42B20; Secondary:
42B30, 42B25, 47B47.\endgraf
{\it Key words and phrases:} product space; Hardy space; Lebesgue space;
sublinear operator; commutator; Calder\'on-Zygmund operator;
Lipschitz function.\endgraf
The
first author is partially supported by a research grant from
United States Army Research Office and a competitive research
grant at Georgetown University.
The second (corresponding) author
is supported by the National
Natural Science Foundation (Grant No. 10871025) of China.}}
\author{\vspace{-0.3cm}
Der-Chen Chang, Dachun Yang\thanks{Corresponding author}\ \ and Yuan Zhou
}
\date{ }
\maketitle

\begin{center}
\begin{minipage}{14cm}
{\small{\bf Abstract.} Let $p\in(0,\,1]$.
In this paper, the authors
prove that a sublinear operator $T$ (which is originally defined on
smooth functions with compact support) can be extended as a
bounded sublinear operator from product Hardy spaces
$H^p({{\mathbb R}^n\times{\mathbb R}^m})$ to
some quasi-Banach space ${\mathcal B}$ if and only if $T$ maps all
$(p,\,2,\,s_1,\,s_2)$-atoms into uniformly bounded elements of
${\mathcal B}$. Here $s_1\ge\lfloor n(1/p-1)\rfloor$
and $s_2\ge\lfloor m(1/p-1)\rfloor$. As usual,
$\lfloor n(1/p-1)\rfloor$ denotes the maximal integer no more than $n(1/p-1)$.
Applying this result, the authors establish the boundedness of the
commutators generated by Calder\'on-Zygmund operators and
Lipschitz functions from the Lebesgue space
$L^p({{\mathbb R}^n\times{\mathbb R}^m})$ with some
$p>1$ or the Hardy space $H^p({{\mathbb R}^n\times{\mathbb R}^m})$
with some $p\le1$ but near 1
to the Lebesgue space $L^q({{\mathbb R}^n\times{\mathbb R}^m})$ with some $q>1$.


}
\end{minipage}
\end{center}


\section{\hspace{-0.6cm}{.} Introduction\label{s1}}

\hskip\parindent The theory of Calder\'on-Zygmund operators and Hardy
spaces on product spaces has been studied by many mathematicians
extensively in the past thirty years, see, for example, \cite{f81,
f85, fe2, fs, hy, j, yz, yz2}. Recently, Ferguson and Lacey
\cite{fl} characterized the product ${\rm BMO}\,(\rr_+^2
\times\rr_+^2)$ by the nested commutator determined by the
one-dimensional Hilbert transform in the $j$th variable, $j=1,2$.
Motivated by this, Chen, Han and Miao in \cite{chm} established the
boundedness on $H^1(\rnm)$ of bi-commutators of fractional integrals
with BMO functions. The boundedness on $H^1(\rnm)$ of the
Marcinkiewicz integral and its commutator with Lipschitz function
was also established in \cite{yz}.

To establish the boundedness of operators on
Hardy spaces on $\rn$ and $\rnm$,
one usually appeals to the atomic decomposition
characterization of Hardy spaces, which means that
a function or distribution in Hardy spaces can be represented as a
linear combination of atoms; see \cite{co,La}
and \cite{cf1,cf3} respectively.
Then, the boundedness of linear operators on
Hardy spaces can be deduced from their behavior on atoms
in principle.

However, Meyer \cite[p.\,513]{mtw} (see also \cite{b,gr}) gave an
example of $f\in H^1(\rn)$, whose norm cannot be achieved by its
finite atomic decompositions via $(1,\,\fz)$-atoms. Based on this
fact, Bownik \cite[Theorem 2]{b} constructed a surprising example
of a linear functional defined on a dense subspace of $H^1(\rn)$,
which maps all $(1,\,\fz)$-atoms into bounded scalars, but yet
cannot extend to a bounded linear functional on the whole
$H^1(\rn)$. This implies that it cannot guarantee the boundedness
of linear operator $T$ from $H^p(\rn)$ with $p\in(0,\,1]$ to some
quasi-Banach space $\cb$ only proving that $T$ maps all
$(p,\,\fz)$-atoms into uniformly bounded elements of  $\cb$. This
phenomenon has also essentially already been observed by Y. Meyer
in \cite[p.\,19]{mc}. Moreover, motivated by this, Yabuta \cite{y}
gave some sufficient conditions for the boundedness of $T$ from
$H^p(\rn)$ with $p\in(0,\,1]$ to $L^q(\rn)$ with $q\ge1$ or
$H^q(\rn)$ with $q\in[p,\,1]$. However, these conditions are not
necessary. In \cite{yz2}, a boundedness criterion was established
as follows: a sublinear operator $T$ (which is originally
defined on smooth functions with compact support) extends to a
bounded sublinear operator from
$H^p(\rn)$ with $p\in(0,1]$ to some quasi-Banach spaces $\cb$ if
and only if $T$ maps all $(p,2)$-atoms into uniformly bounded
elements of $\cb$. This result shows the structure difference
between atomic characterization of $H^p(\rn)$ via $(p,2)$-atoms
and $(p,\fz)$-atoms. This result is generalized to spaces of
homogeneous type in \cite{yz3}.

The purpose of this paper is two folds. We first generalize the
boundedness criterion on $\rn$ in \cite{yz2} to product Hardy
spaces on $\rnm$. Precisely, we prove that a sublinear operator $T$
(which is originally defined on smooth functions
with compact support) extends to a bounded sublinear operator
from $H^p(\rnm)$ with $p\in(0,1]$ to some
quasi-Banach spaces $\cb$ if and only if $T$ maps all
$(p,2)$-atoms into uniformly bounded elements of $\cb$. Invoking
this result and motivated by \cite{chm,fl,yz}, we then establish
the boundedness of the commutators generated by Calder\'on-Zygmund
operators and Lipschitz functions from the Lebesgue space
$L^p(\rnm)$ with some $p>1$ or the Hardy space $H^p(\rnm)$ with
some $p\le1$ but near $1$ to the Lebesgue space $L^q(\rnm)$ with
some $q>1$.

To state the main results, we first recall some notation and
notions on product Hardy spaces. For $n,\,m\in\nn$, denote
by $\cs(\rnm)$ the space of Schwartz functions on $\rnm$ and
by $\cs'(\rnm)$ its dual space. Let $\cd(\rnm)$ be the space of
all smooth functions on $\rnm$ with compact support. For
$s_1,\,s_2\in\zz_+$, let $\cd_{s_1,\,s_2}(\rnm)$ be the set of all
functions $f\in\cd(\rnm)$ with vanishing moments up to order $s_1$
with respect to the first variable and order $s_2$ with respect to
the second variable. More precisely, if $f\in\cd(\rnm)$, then for
$\az_1\in\zz_+^{n}$ and
$\az_2\in\zz_+^{m}$ with $|\az_1|\le s_1$ and $|\az_2|\le s_2$,
one has
$$\begin{array}{cl}
\dint_\rn f(x_1,\,x_2)x_1^{\az_1}
\,dx_1&=0\quad{\mbox{for
all}}\quad x_2\in\rrm,\\
\dint_\rrm f(x_1,\,x_2)x_2^{\az_2}\,dx_2&=0\quad{\mbox{for
all}}\quad x_1\in\rn.
\end{array}
$$
For $s_1,\,s_2\in\zz_+$ and $\sz_1,\,\sz_2\in[0,\fz)$, we denote by
$\cd_{s_1,\,s_2;\,\sz_1,\,\sz_2}(\rnm)$ the space
$\cd_{s_1,\,s_2}(\rn\times\rrm)$ endowed with the norm
$$\|f\|_{\cd_{s_1,\,s_2;\,\sz_1,\,\sz_2}(\rnm)}
\equiv\sup_{x_1\in\rn,\,x_2\in\rrm}(1+|x_1|)^{\sz_1}(1+|x_2|)^{\sz_2}|f(x_1,\,x_2)|.$$

In articles \cite{cf1,cf2,cf3}, Chang and Fefferman introduced the
following atoms and atomic Hardy spaces on the product space
$\rnm$.

\medskip

\begin{defn}\hspace{-0.2cm}{.}\rm\label{d1.1}
Let $p\in(0,1]$, $s_1\ge\laz n(1/p-1)\raz $ and $s_2\ge \laz m(1/p-1)\raz $. A
function $a$ supported in an open set $\boz\subset\rnm$ with
finite measure is said to be a $(p,\,2,\,s_1,\,s_2)$-atom provided
that
\begin{enumerate}
\item[(AI)] $a$ can be written as $a=\sum_{R\in \cm(\boz)}a_R,$
where $\cm(\boz)$ denotes all the maximal dyadic subrectangles of
$\boz$ and $a_R$ is a function satisfying that
\begin{enumerate}
\item[(i)] $a_R$ is supported on $2R=2I\times 2J$, which is a
rectangle with the same center as $R$ and whose side length is $2$
times that of R, \item[(ii)] $a_R$ satisfies the cancelation
conditions that
$$\begin{array}{cl}
\dint_{2I}a_R(x_1,\,x_2)x_1^{\az_1}\,dx_1=0&\quad
{for\ all}\ x_2\in 2J\ and\ |\az_1|\le s_1,\\
\dint_{2J}a_R(x_1,\,x_2)x_2^{\az_2}\,dx_2=0 &\quad {for\ all}\
x_1\in 2I\ and\ |\az_2|\le s_2;
\end{array}$$
\end{enumerate}
\item[(AII)] $a$ satisfies the size conditions that
$\|a\|_{L^2(\rnm)}\le |\boz|^{1/2-1/p}$ and
$$\lf(\sum_{R\in\cm(\boz)}\|a_R\|^2_{L^2(\rnm)}\r)^{1/2}\le
|\boz|^{1/2-1/p}.$$
\end{enumerate}
\end{defn}

\begin{defn}\hspace{-0.2cm}{.}\rm\label{d1.2}
Let $p\in(0,\,1]$, $s_1\ge\laz n(1/p-1)\raz $ and $s_2\ge \laz
m(1/p-1)\raz $. A distribution $f\in\cs'(\rnm)$ is said to be an
element in $H^{p,\,2,\,s_1,\,s_2}(\rnm)$ if there exist a sequence
$\{\lz_k\}_{k\in\nn}\subset\cc$ and $(p,\,2,\,s_1,\,s_2)$-atoms
$\{a_k\}_{k\in\nn}$ such that $f=\sum_{k\in\nn}\lz_k a_k$ in
$\cs'(\rnm)$ with $\sum_{k\in\nn}|\lz_k|^p<\fz$. Moreover, define
the quasi-norm of $f\in H^{p,\,2,\,s_1,\,s_2}(\rnm)$  by
$\|f\|_{H^{p,\,2,\,s_1,\,s_2}(\rnm)}\equiv \inf \{
(\sum_{k\in\nn}|\lz_k|^p )^{1/p} \},$ where the infimum is taken
over all the decompositions as above.
\end{defn}

It is well known that
$H^{p,\,2,\,s_1,\,s_2}(\rnm)=H^{p,\,2,\,t_1,\,t_2}(\rnm)$ with
equivalent norms when $s_1,\,t_1\ge\laz n(1/p-1)\raz $ and
$s_2,\,t_2\ge\laz m(1/p-1)\raz $; see \cite{cf1,cf2,cf3,fe1,h89}. Thus, we
denote $H^{p,\,2,\,s_1,\,s_2}(\rnm)$ simply by $H^p(\rnm)$.

Recall that a quasi-Banach space $\cb$ is a vector space endowed with
a quasi-norm $\|\cdot\|_\cb$ which is nonnegative, non-degenerate
(i.\,e., $\|f\|_\cb=0$ if and only if $f=0$), homogeneous,
and obeys the quasi-triangle inequality, i.\,e., there exists
a constant $C_0\ge1$ such that for all $f,\, g\in\cb$,
$$\|f+g\|_\cb\le C_0(\|f\|_\cb+\|g\|_\cb).\leqno(1.1)$$

\begin{defn}\hspace{-0.2cm}{.}\rm\label{d1.3}
Let $q\in(0,\,1]$. A quasi-Banach spaces $\cbq$ with
the quasi-norm $\|\cdot\|_\cbq$ is said to be a
$q$-quasi-Banach space if $\|\cdot\|_\cbq^q$
satisfies the triangle inequality, i.\,e.,
$\|f+g\|^q_\cbq\le \|f\|^q_\cbq+\|g\|^q_\cbq$
for all $f,\, g\in\cbq$.
\end{defn}

We point out that by the Aoki theorem (see \cite{a} or
\cite[p.\,66]{g}), any quasi-Banach space with the positive constant
$C_0$ as in (1.1) is essentially a $q$-quasi-Banach space with
$q=\lfloor\log_2(2C_0)\rfloor^{-1}$. From this, any Banach space is
a $1$-quasi-Banach space. Moreover, $\ell^q$, $L^q(\rnm)$ and
$H^q(\rnm)$ with $q\in(0,\,1)$ are typical $q$-quasi-Banach spaces.

Let $q\in(0,\,1]$. For any given $q$-quasi-Banach space $\cbq$  and
linear space ${\mathcal Y}$, an operator $T$ from ${\mathcal Y}$
to $\cbq$ is called to be $\cbq$-sublinear if for any $f,\,
g\in{\mathcal Y}$ and $\lz,\,\nu\in\cc$, we have
$$\|T(\lz f+\nu g)\|_\cbq
\le\lf(|\lz|^q\|T(f)\|^q_\cbq+|\nu|^q\|T(g)\|^q_\cbq\r)^{1/q}$$
and $\|T(f)-T(g)\|_\cbq\le \|T(f-g)\|_\cbq$; see \cite{yz2, yz3}.
Obviously, if $T$ is
linear, then $T$ is $\cbq$-sublinear. Moreover, if $\cbq$
is a space of functions, $T$ is sublinear in the classical
sense and $T(f)\ge 0$ for all $f\in {\mathcal Y}$,
then $T$ is also $\cbq$-sublinear.

The following is one of main results in this paper,
which generalizes the main result in \cite{yz2} to product Hardy spaces.

\begin{thm}\hspace{-0.2cm}{.}\label{t1.1}
Let $p\in(0,\,1]$, $q\in[p,\,1]$ and $\cbq$ be a $q$-quasi-Banach
space. Suppose that $s_1\ge\laz n(1/p-1)\raz $ and $s_2\ge\laz m(1/p-1)\raz $. Let
$T$ be a $\cb_q$-sublinear operator from $\cd_{s_1,\,s_2}(\rnm)$ to $\cbq$.
Then $T$ can be extended as a bounded $\cb_q$-sublinear operator from $H^p(\rnm)$ to
$\cbq$ if and only if
 $T$ maps all  $(p,\,2,\,s_1,\,s_2)$-atoms in $\cd_{s_1,\,s_2}(\rnm)$
 into uniformly bounded elements of $\cbq$.
\end{thm}

Theorem 1.1 further complements the proofs of Theorem 1 in
\cite{fe2} and a theorem in \cite{f85}, whose proof is
presented in Section 2 below. The necessity of Theorem 1.1 is
obvious. To prove the sufficiency, for $p\in(0,1]$,
$s_1\ge\laz n(1/p-1)\raz$, $s_2\ge \laz m(1/p-1)\raz $ and
$f\in\cd_{s_1,\,s_2}(\rnm)$, we first prove that $f$ has an atomic
decomposition which converges in
$\cd_{s_1,\,s_2;\,\sz_1,\,\sz_2}(\rnm)$ for some
$\sz_1\in(\max\{n/p,\,n+s\},\,n+s+1)$ and
$\sz_2\in(\max\{n/p,\,n+s\},\,n+s+1)$ (Lemma 2.3), and then extend
$T$ to the whole $\cd_{s_1,\,s_2;\,\sz_1,\,\sz_2}(\rnm)$ boundedly
(Lemma 2.4). Finally, we continuously extend $T$ to the whole
$H^p(\rnm)$ by using the density of
$\cd_{s_1,\,s_2}(\rnm)$ in $H^p(\rnm)$.

Recall that a function $a$ is said to be a rectangular
$(p,\,2,\,s_1,\,s_2)$-atom if
\begin{enumerate}
\item[(R1)] $\supp a \subset R=I\times J$, where $I$ and $J$ are
cubes in $\rn$ and $\rrm$, respectively; \item[(R2)] $\int_\rn
a(x_1,\,x_2)x_1^{\az_1}\,dx_1=0$ for all $x_2\in\rrm$ and
$|\az_1|\le s_1$, and $\int_\rrm a(x_1,\,x_2)x_2^{\az_2}\,dx_2=0$
for all $x_1\in\rn$ and $|\az_2|\le s_2$; \item[(R3)]
$\|a\|_{L^2(\rnm)}\le|R|^{1/2-1/p}$.
\end{enumerate}

As a consequence of Theorem 1.1, we obtain the following result
which includes a fractional version of Theorem 1 in \cite{fe2}
and is known to have many applications in harmonic analysis.

\begin{cor}\hspace{-0.2cm}{.}\label{c1.1}
Let $q_0\in[2,\fz)$ and $T$ be a bounded sublinear operator from
$L^2(\rnm)$ to $L^{q_0}(\rnm)$. Let $p\in(0,1]$ and
$1/q-1/p=1/q_0-1/2$. If there exist positive constants $C$ and $\dz$
such that for all rectangular $(p,\,2,\,s_1,\,s_2)$-atoms $a$
supported in $R$ and all $\gz\ge 8\max\{n^{1/2},\,m^{1/2}\}$,
$$\int_{(\rnm)\setminus\wz R_\gz} |Ta(x_1,\,x_2)|^q\,dx_1\,dx_2\le C\gz^{-\dz},$$
where $\wz R_\gz$ denotes the $\gz$-fold enlargement of $R$, then
$T$ can be extended as a bounded sublinear operator from $H^p(\rnm)$ to
$L^q(\rnm)$.
\end{cor}

The proof of Corollary 1.1 is given in Section 2 below.
We point out that if $q_0=2$ and $T$ is linear, then Corollary 1.1
is just Theorem 1 in \cite{fe2}. Moreover, there exists a gap
in the proof of Theorem 1 in \cite{fe2} (so is the proof
of a theorem in \cite{f85}), namely, it was not clear in \cite{fe2}
how to deduce the boundedness of the considered linear operator
$T$ on the whole Hardy space $H^p(\rnm)$ from its boundedness
uniformly on atoms. Our Theorem 1.1 here seals this gap.

\begin{rem}\hspace{-0.2cm}{.}\rm\label{r1.1}
Using Corollary 1.1, we now give affirmative answers to
the questions in Remark 4.2 and
Remark 4.3 of \cite{yz}. We use the same notation and
notions as in \cite{yz}. Particularly,
denote by $\mu_\boz$ the Marcinkiewicz integral operator on $\rnm$
with kernel $\boz\in\lip(\az_1,\az_2; {\mathbb S}^{n-1},{\mathbb S}^{m-1})$,
here $\az_1,\az_2\in(0,1]$. If
$\max\{n/(n + \az_1), m/(m + \az_2)\}<p\le1$, then
in Remark 4.2 of \cite{yz}, we proved that for all
$(p,2,0,0)$ atoms $a$, $\|\mu_\boz(a)\|_{L^p(\rnm)}\ls1$.
Moreover, let $b\in\lip(\bz_1,\bz_2;\rnm)$ with
$\bz_1, \bz_2\in(0,1]$ satisfying $\bz_1/n =\bz_2/m$ and
$C_b(\mu_\boz)$ be the commutator of $b$ and $\mu_\boz$.
If $1/q=1/p-\bz_1/n$ and
$$\max\{n/(n + \az_1), m/(m + \az_2)\}<p\le1,$$
then in Remark 4.3 of \cite{yz}, we proved that for all
$(p,2,0,0)$ atoms $a$,
$$\|C_b(\mu_\boz)(a)\|_{L^q(\rnm)}\ls1.$$
However, in \cite{yz}, it is not clear how to obtain the boundedness
of  $\mu_\boz$ from $H^p(\rnm)$
to  $L^p(\rnm)$ and boundedness
of  $C_b(\mu_\boz)$ from $H^p(\rnm)$
to  $L^q(\rnm)$ by these known facts.
Applying Theorem 1.1 here, we now obtain
these desired boundedness, and hence answer the questions
in Remark 4.2 and Remark 4.3 of \cite{yz}.
\end{rem}

Now we turn to the boundedness of commutators generated by Lipschitz
functions and Calder\'on-Zygmund operators. We first introduce the
notion of Lipschitz functions on $\rnm$. Let $\az\in (0, 1]$. A function
$b$ on $\rn$ is said to belong to $\lip(\az;\,\rn)$ if there exists
a positive constant $C$ such that for all $x,\,x'\in\rn$,
$$|b(x)-b(x')|\le C|x-x'|^\az.$$
Obviously, a function in the space $\lip(\az;\,\rn)$ is not necessary
bounded. For example, $|x|^\az\in\lip(\az;\,\rn)$, but
$|x|^\az\not\in L^\fz(\rn)$.

\begin{defn}\hspace{-0.2cm}{.}\rm\label{d1.4}
Let $\az_1,$ $\az_2\in (0, 1]$. A function $f$ on $\rnm$ is
said to belong to $\lip(\az_1,\,\az_2;\,\rnm)$, if there exists a
positive constant $C$ such that for all $x_1,\, y_1\in\rn$ and
$x_2,\, y_2\in\rr^m$,
$$|[f(x_1,\,x_2)-f(x_1,\,y_2)]-[f(y_1,\,x_2)-f(y_1,\,y_2)]|
\le C|x_1-y_1|^{\az_1}|x_2-y_2|^{\az_2}.\leqno(1.2)$$ The minimal
constant $C$ satisfying (1.2) is defined to be the norm of $f$ in
the space $\laa$ and denoted by $\|f\|_\laa$.
\end{defn}

We remark that a function in the space $\lip(\az_1,\,\az_2;\,\rnm)$
is also not necessary to be bounded. In fact, if $f_1\in
\lip(\az_1;\,\rn)$ and $f_2\in\lip(\az_2;\,\rr^m)$, then it is easy
to check $f_1(x_1)f_2(x_2)\in\lip(\az_1,\,\az_2;\,\rnm)$.

In this paper, we consider a class of Calder\'on-Zygmund operators
$T$ on $\rnm$, whose kernel $K$ is a continuous function on
$(\rn\times\rn\times\rr^m\times\rr^m)
\setminus\{(x_1,\,y_1,\,x_2,\,y_2):\  x_1=y_1\ or \ x_2=y_2\}$ and
satisfies that there exist positive constants $C$ and
$\ez_1,\,\ez_2\in (0, 1]$ such that
\begin{enumerate}
\item[(K1)] for all $x_1\ne y_1$ and $x_2\ne y_2$,
$$|K(x_1,\,y_1,\,x_2,\,y_2)|\le C\frac
1{|x_1-y_1|^n}\frac 1{|x_2-y_2|^m};$$
\item[(K2)] for all
$x_1\ne y_1$, $x_2\ne y_2$, $z_1\in\rn$ and
$|y_1-z_1|\le|x_1- y_1|/2$,
$$ |K(x_1,\,y_1,\,x_2,\,y_2)-K(x_1,\,z_1,\,x_2,\,y_2)|\le
C\dfrac{|y_1-z_1|^{\ez_1}}{|x_1-y_1|^{n+\ez_1}}\dfrac{1}{|x_2-y_2|^{m}};$$
\item[(K3)]
for all $x_1\ne y_1$, $x_2\ne y_2$, $z_2\in\rr^m$ and $|y_2-z_2|\le
|x_2- y_2|/2$,
$$ |K(x_1,\,y_1,\,x_2,\,y_2)-K(x_1,\,y_1,\,x_2,\,z_2)|\le
C\dfrac{1}{|x_1-y_1|^n}\dfrac{|y_2-z_2|^{\ez_2}}{|x_2-y_2|^{m+\ez_2}};$$
\item[(K4)] for all $x_1\ne y_1$,\ $x_2\ne y_2$,\ $z_1\in\rn$, $z_2\in\rr^m$,
$|y_1-z_1|\le|x_1- y_1|/2$ and $|y_2-z_2|\le|x_2- y_2|/2$,
$$\begin{array}{l}
\lf|[K(x_1,\,y_1,\,x_2,\,y_2)-K(x_1,\,z_1,\,x_2,\,y_2)]
-[K(x_1,\,y_1,\,x_2,\,z_2)-K(x_1,\,z_1,\,x_2,\,z_2)]\r|\\
\quad\le
C\dfrac{|y_1-z_1|^{\ez_1}}{|x_1-y_1|^{n+\ez_1}}
\dfrac{|y_2-z_2|^{\ez_2}}{|x_2-y_2|^{m+\ez_2}}.
\end{array}$$
\end{enumerate}
The minimal constant $C$ satisfying (K1) through (K4) is denoted by $\|K\|$.

Let $\az_1,\, \az_2\in (0, 1]$,\ $b\in \lip(\az_1,\,\az_2;\,\rnm)$ and
$T$ be any Calder\'on-Zygmund operator with kernel $K$ satisfying
the above conditions from (K1) to (K4). For any suitable function
$f$ and $(x_1,x_2)\in\rnm$, define the commutator $[b,\,T]$ by
$$\begin{array}[t]{cl}
[b,\,T](f)(x_1,\,x_2)
&=\dint_{\rnm} K(x_1,\,y_1,\,x_2,\,y_2)
\\
&\quad\times[b(x_1,\,x_2)-b(x_1,\,y_2)-b(y_1,\,x_2)+b(y_1,\,y_2)]f(y_1,\,y_2)\,
d y_1 d y_2.\end{array}\leqno(1.3)$$

The following result gives the boundedness of the commutator $[b,\,T]$
on Lebesgue spaces.

\begin{thm}\hspace{-0.2cm}{.}\label{t1.2}
Let $\ez_1,\,\ez_2,\,\az_1,\,\az_2\in (0, 1]$, $\az_1/n=\az_2/m$,
$p\in (1, n/\az_1)$ and $1/q=1/p-\az_1/n$. Let
$b\in\lip(\az_1,\,\az_2;\,\rnm)$, $T$ be a Calder\'on-Zygmund
operator whose kernel $K$ satisfies the conditions from (K1) to
(K4), and $[b,\,T]$ be the commutator as in (1.3). Then there exists
a positive constant $C $ independent of $\|b\|_\laa$ and $\|K\|$
such that for all $f\in \lp$,
$$\|[b,\,T] (f)\|_\lq\le C\|K\|\|b\|_\laa\|f\|_\lp.$$
\end{thm}

Here is another main result of this paper, whose proof depends on
Corollary 1.1.

\begin{thm}\hspace{-0.2cm}{.}\label{t1.3}
Let $0<\az_1\le\min\{n/2,\,1\}$, $\az_1/n=\az_2/m$,
$\ez_1,\,\ez_2\in (0, 1]$,
$$\max\{n/(n+\ez_1),\,n/(n+\az_1),\,m/(m+\ez_2),
\,m/(m+\az_2)\}<p\le1\leqno(1.4)$$ and $1/q=1/p-\az_1/n.$ Assume that
$b\in\laa$. Let $T$ be a Calder\'on-Zygmund operator whose kernel
$K$ satisfies the conditions (K1) through (K4), and $[b,\, T]$ be
the commutator defined in (1.3). Then there exists a positive
constant $C$ independent of $\|b\|_\laa$ and $\|K\|$ such that for
all $f\in \hp$,
$$\|[b,\,T] (f)\|_\lq\le C\|K\|\|b\|_\laa\|f\|_\hp.$$
\end{thm}

The proofs of Theorem 1.2 and Theorem 1.3 are presented in Section 3.

We finally make some conventions. Throughout this paper,
let $\nn=\{1,\,2,\,\cdots\}$ and $\zz_+=\nn\cup\{0\}$.
We always use $C$ to denote a positive constant that is independent of main
parameters involved but whose value may differ from line to line.
We use $f\ls g$ to denote $f\le Cg$ and $f\sim g$ to denote $f\ls g\ls f$.

\section{\hspace{-0.6cm}{\bf .} Proofs of Theorem 1.1 and Corollary 1.1}\label{s2}

\hskip\parindent As a matter of convenience, in this section, we denote $n$ and $m$,
respectively, by $n_1$ and $n_2$. For $i=1,\,2$ and
$s_i\in\zz_+$, denote by $\cd_{s_i}(\rr^{n_i})$ the set of all
smooth functions with compact support and vanishing moments up to
order $s_i$. Then there exist functions
$\psi^{(i)}\in\cd_{s_i}(\rr^{n_i})$ and $\varphi^{(i)}\in\cs(\rr^{n_i})$ such that

(i) $\supp \psi^{(i)}\subset B^{(i)}(0,\,1)$, $\widehat {\psi^{(i)}}\ge 0$
and $\widehat {\psi^{(i)}}(\xi_i)\ge \frac 12$ if $\frac 12\le|\xi_i|\le 2$,
where and in what follows $B^{(i)}(0,\,r_i)\equiv\{x_i\in\rr^{n_i}:\
|x_i|<r_i\}$ and $\widehat{\psi^{(i)}}$ denotes the Fourier
transform of $\psi^{(i)}$;

(ii) $\supp\widehat {\vz^{(i)}}\subset \{\xi_i\in\rr^{n_i}:\ 1/2\le|\xi_i|\le
2\}$ and $\widehat {\vz^{(i)}}\ge 0$;

(iii) $\sup\{\widehat {\vz^{(i)}}(\xi_i):\ 3/5\le|\xi_i|\le 5/3\}>C$
for some positive constant $C$;

(iv)
$\int_0^\fz\widehat {\psi^{(i)}}(t_i\xi_i)\widehat {\vz^{(i)}}(t_i\xi_i)
\,\frac{dt_i}{t_i}=1$
for all $\xi_i\in{\rr^{n_i}}\setminus\{0\}$.

 \noindent Such $\psi^{(i)}$ and
$\vz^{(i)}$ can be constructed by a slight modification of Lemma
(1.2) of \cite{fjw}; see also Lemma (5.12) in \cite{fjw} for a
discrete variant. Then by an argument similar to the proofs of
Theorem (1.3) and  Theorem 1 in Appendix of \cite{fjw}, we have that
for all $f\in\cs(\rnt)$ and $(x_1,x_2)\in\rnt$,
$$f(x_1,\,x_2)=\int_0^\fz\int_0^\fz(\psi_{t_1,\,t_2}\ast
\varphi_{t_1,\,t_2}\ast
f)(x_1,\,x_2)\,\frac{dt_1}{t_1}\frac{dt_2}{t_2}\leqno(2.1)$$
in both $L^2(\rnt)$ and pointwise, where
and in what follows, for any $i=1,\,2$, $\phi^{(i)}\in\cs(\rr^{n_i})$,
 $x_i\in\rr^{n_i}$ and
$t_i\in(0,\,\fz)$, we always let $\phi^{(i)}_{t_i}(x_i)\equiv
t_i^{-n_i}\phi^{(i)}(t_i^{-1}x_i)$ and
$\phi_{t_1,\,t_2}(x_1,x_2)\equiv\phi^{(1)}_{t_1}(x_1)\phi^{(2)}_{t_2}(x_2)$.
For any set $E\subset(\rnm)$, set $E^\complement\equiv(\rnm)\setminus E$.

\begin{lem}\hspace{-0.2cm}{.}\label{l2.1}
Let $s_i\in\zz_+$, $\psi^{(i)}\in\cd_{s_i}(\rr^{n_i})$ and
$\vz^{(i)}\in\cs(\rr^{n_i})$ satisfy the above conditions (i)
through (iv), where $i=1,\,2$.
Let $0<\sz_i<\sz_i'<n_i+s_i+1$ for $i=1,\,2$.
Then for any $f\in\cd_{s_1,\,s_2}(\rnt)$, there exists a positive
constant $C$ such that
for all $\ez_1,\,\ez_2\in(0,\,1)$ and $L_1, L_2\in(1,\fz)$,
\begin{eqnarray*}
&&\sup_{(x_1,\,x_2)\in\rnt}(1+|x_1|)^{\sz_1}(1+|x_2|)^{\sz_2}\\
&&\quad\quad\times\lf(\int_0^{\ez_1}\int_0^\fz
+\int_{L_1}^\fz\int_0^\fz+ \int_0^\fz\int_0^{\ez_2}
+\int_0^\fz\int_{L_2}^\fz\r)\int_\rnt|(\varphi_{t_1,\,t_2}\ast f)(y_1,\,y_2)|\\
&&\quad\quad\times |\psi_{t_1,\,t_2}(x_1-y_1,\,x_2-y_2)|
\,dy_1\,dy_2\frac{dt_1}{t_1}\frac{dt_2}{t_2}\\
&&\quad\le
C\lf[\ez_1+\ez_2+(L_1)^{\sz_1-n_1-s_1-1}+(L_2)^{\sz_2-n_2-s_2-1}\r],
\end{eqnarray*}
$$\begin{array}[t]{ccl}
&&\dsup_{(x_1,\,x_2)\in\rnt}(1+|x_1|)^{\sz_1}(1+|x_2|)^{\sz_2}\\
&&\quad\quad\times\dint_0^{L_1}\int_0^\fz
\dint_{[B^{(1)}(0,\,2L_1)]^\complement\times\rr^{n_2}}
|(\varphi_{t_1,\,t_2}\ast f)(y_1,\,y_2)|\\
&&\quad\quad\times
 |\psi_{t_1,\,t_2}(x_1-y_1,\,x_2-y_2)|
\,dy_1\,dy_2\frac{dt_1}{t_1}\frac{dt_2}{t_2}
\le C(L_1)^{\sz_1-\sz_1'}
\end{array}\leqno(2.2)$$
and (2.2) with $L_1$, $\sz_1$, $n_1$, $s_1$ and $B^{(1)}$ replaced,
respectively, by $L_2$, $\sz_2$, $n_2$, $s_2$ and $B^{(2)}$.
\end{lem}

In order to prove Lemma 2.1, we need the following technical lemma.
For $i=1,\,2$, $u_i\ge0$, let
$$\cs_{u_i}(\rr^{n_i})\equiv\lf\{\vz\in\cs(\rr^{n_i}):\
\int_{\rr^{n_i}}\vz(x_i)x_i^\az\,dx_i=0,\ |\az|\le u_i\r\}.$$
 For any
$s_1,\,s_2\in\zz_{-1}\equiv\nn\cup \{0,-1\}$, we denote by
$\cs_{s_1,\,s_2}(\rnt)$ the space of functions in $\cs(\rnt)$ with
the vanishing moments up to order $s_1$ in the first variable and
order $s_2$ in the second variable, where we say that
$f\in\cs(\rnt)$ has vanishing moments up to order $-1$ in the
first or second variable, if $f$ has no vanishing moment with
respect to that variable.

\begin{lem}\hspace{-0.2cm}{.}\label{l2.2}
Let $s_i\in\zz_{-1}$, $u_i\in\zz_{-1}$, $\sz_i\in[0,\,\fz)$ and
$\vz^{(i)}\in\cs_{u_i}(\rr^{n_i})$ for $i=1,\,2$. For any
$f\in\cs_{s_1,\,s_2}(\rnt)$, there exists a positive constant $C$ such
that
\begin{enumerate}
\item[(i)] if $u_1>-1$, then for all $t_1\in(0,\, 1]$ and $(x_1,\,x_2)\in\rnt$,
$$|(\varphi^{(1)}_{t_1}\ast_1f)(x_1,\,x_2)|\le C
t_1^{u_1+1}(1+|x_1|)^{-\sz_1} (1+|x_2|)^{-\sz_2},$$
where and in what follows
$(\varphi^{(1)}_{t_1}\ast_1f)(x_1,\,x_2)\equiv\int_{\rr^{n_1}}
\varphi^{(1)}_{t_1}(y_1) f(x_1-y_1,\,x_2)\,dy_1;$
\item[(ii)] if
$s_1>-1$, then for all $t_1\in[1,\, \fz)$ and
$(x_1,\,x_2)\in\rnt$,
$$|(\varphi^{(1)}_{t_1}\ast_1f)(x_1,\,x_2)|
\le
Ct_1^{-n_1-s_1-1}\lf(1+\frac{|x_1|}{t_1}\r)^{-\sz_1}(1+|x_2|)^{-\sz_2};$$
\item[(iii)] if $u_1,\,u_2>-1$, then for all $t_1,\,t_2\in(0, 1]$ and
$(x_1,\,x_2)\in\rnt$,
$$|(\varphi_{t_1,\,t_2}\ast f)(x_1,\,x_2)|
\le C t_1^{u_1+1}t_2^{u_2+1}(1+|x_1|)^{-\sz_1} (1+|x_2|)^{-\sz_2};$$ \item[(iv)]
if $u_1,\,s_2>-1$, then for all $t_1\in(0,\,1]$, $t_2\in[1,\,\fz)$ and
$(x_1,\,x_2)\in\rnt$,
$$|(\varphi_{t_1,\,t_2}\ast f)(x_1,\,x_2)|\le
Ct_1^{u_1+1}t_2^{-n_2-s_2-1}(1+|x_1|)^{-\sz_1}\lf(1+\frac{|x_2|}{t_2}\r)^{-\sz_2};$$
\item[(v)] if $s_1,\,u_2>-1$, then for all $t_1\in[1,\,\fz)$,
$t_2\in(0,\,1]$ and $(x_1,\,x_2)\in\rnt$,
$$|(\varphi_{t_1,\,t_2}\ast f)(x_1,\,x_2)|\le
Ct_1^{-n_2-s_2-1}t_2^{u_2+1}\lf(1+\frac{|x_1|}{t_1}\r)^{-\sz_1}(1+|x_2|)^{-\sz_2};$$
\item[(vi)] if $s_1,\,s_2>-1$, then for all
$t_1,\,t_2\in[1,\,\fz)$ and $(x_1,\,x_2)\in\rnt$,
$$|(\varphi_{t_1,\,t_2}\ast f)(x_1,\,x_2)|\le
Ct_1^{-n_1-s_1-1}t_2^{-n_2-s_2-1}\lf(1+\frac{|x_1|}{t_1}\r)^{-\sz_1}
\lf(1+\frac{|x_2|}{t_2}\r)^{-\sz_2}.$$
\end{enumerate}
\end{lem}

\begin{pf}\rm
To prove Lemma 2.2, we use some ideas in the proofs of Lemma 2 and
Lemma 4 in Appendix (III) of \cite{fjw}.

To prove (i), by $\int_{\rr^{n_1}} \varphi^{(1)}(x_1)x_1^\az\,dx_1=0$
for $|\az|\le u_1$, we
have
\begin{eqnarray*}
(\varphi_{t_1}\ast_1 f)(x_1,\,x_2)
&&=\int_{\rr^{n_1}}\vz^{(1)}_{t_1}(y_1)\lf[f(x_1-y_1,\,x_2)
-\sum_{|\gz|\le u_1}\frac{1}{\gz!}y_1^\gz
(D_1^\gz f)(x_1,\,x_2)\r]dy_1\\
&&=\int_{|y_1|<|x_1|/2}\vz^{(1)}_{t_1}(y_1)\lf[f(x_1-y_1,\,x_2)
-\sum_{|\gz|\le u_1}\frac{1}{\gz!}y_1^\gz
(D_1^\gz f)(x_1,\,x_2)\r]dy_1\\
&&\hs+\int_{|y_1|\ge|x_1|/2}\cdots\\
&&\equiv I_1+I_2.
\end{eqnarray*}
For the estimation of $I_1$, noticing that $|x_1|/2\le
|x_1-z_1|\le2|x_1|$ for $|z_1|\le |x_1|/2$, by $|y_1|<|x_1|/2$ and
the mean value theorem, we obtain
$$\begin{array}[t]{ccl}
&&\lf|f(x_1-y_1,\,x_2)-\dsum_{|\gz|\le u_1}\frac{1}{\gz!}y_1^\gz
(D_1^\gz f)(x_1,\,x_2)\r|\\
&&\quad=\dsup_{|\gz|={u_1+1}}\dsup_{|z_1|\le
|x_1-y_1|}|(D^\gz_1 f)(x_1-z_1,\,x_2)|
|y_1|^{u_1+1}\\
&&\quad\ls |y_1|^{u_1+1}\dsup_{|z_1|\le|x_1|/2}(1+|x_1-z_1|)^{-\sz_1}(1+|x_2|)^{-\sz_2}\\
&&\quad\ls|y_1|^{u_1+1}(1+|x_1|)^{-\sz_1}(1+|x_2|)^{-\sz_2},
\end{array}\leqno(2.3)$$
where $\gz=(\gz_1,\,\cdots,\,\gz_{n_1})\in\zz_+^{n_1}$,
$x_1=(x^1_1,\,\cdots,\,x^{n_1}_1)$
and $D^\gz_1=(\frac{\partial}{\partial x^1_1})^{\gz_1}
\cdots(\frac{\partial}{\partial x^{n_1}_1})^{\gz_{n_1}}$.
This leads to that
\begin{eqnarray*}
|I_1|&&\ls (1+|x_1|)^{-\sz_1}(1+|x_2|)^{-\sz_2}
\int_{|y_1|<|x_1|/2}|y_1|^{u_1+1}|\varphi_{t_1}^{(1)}(y_1)|\,dy_1\\
&&\ls t_1^{u_1+1}(1+|x_1|)^{-\sz_1}(1+|x_2|)^{-\sz_2}
\int_{\rr^{n_1}}|y_1|^{u_1+1}|\varphi^{(1)}(y_1)|\,dy_1\\
&&\ls t_1^{u_1+1}(1+|x_1|)^{-\sz_1}(1+|x_2|)^{-\sz_2}.
\end{eqnarray*}
To estimate $I_2$, similarly to (2.3), we have
$$\begin{array}[t]{ccl}
&&\lf|f(x_1-y_1,\,x_2)-\dsum_{|\gz|\le s_1}\frac{1}{\gz!}y_1^\gz
(D_1^\gz f)(x_1,\,x_2)\r|\ls|y_1|^{u_1+1}(1+|x_2|)^{-\sz_2}.
\end{array}\leqno(2.4)$$
If $|x_1|\ge1$ and $\sz_1>0$,
by $|x_1|^{-1}\le 2(1+|x_1|)^{-1}$ and (2.4), for all
$t_1\in(0,\,1]$, we have
\begin{eqnarray*}
|I_2|&&\ls(1+|x_2|)^{-\sz_2}
\int_{|y_1|\ge|x_1|/2}|y_1|^{u_1+1}|\varphi^{(1)}_{t_1}(y_1)|\,dy_1\\
&&\ls(1+|x_2|)^{-\sz_2}
t_1^{u_1+1}\int_{|y_1|\ge|x_1|/(2t_1)}|y_1|^{u_1+1}|\varphi^{(1)}(y_1)|\,dy_1\\
&&\ls t_1^{u_1+1}(1+|x_2|)^{-\sz_2}\int_{|x_1|/(2t_1)}^\fz
r_1^{-\sz_1-1}\,dr_1\\
&&\ls t_1^{u_1+1}(1+|x_1|)^{-\sz_1}(1+|x_2|)^{-\sz_2}.
\end{eqnarray*}
If $|x_1|\le1$ or $\sz_1=0$, by (2.4),
$$|I_2|\ls t_1^{u_1+1}(1+|x_2|)^{-\sz_2}\int_{\rr^{n_1}}
|y_1|^{u_1+1}|\varphi^{(1)}(y_1)|\,dy_1\ls
t_1^{u_1+1}(1+|x_1|)^{-\sz_1}(1+|x_2|)^{-\sz_2}.$$
Thus combining the
estimations for $I_1$ and $I_2$ yields (i).

To prove (ii), since $\varphi^{(1)}\in\cs_0(\rr^{n_1})$ and
$f\in\cs_{s_1,\,s_2}(\rnt)$, we have
\begin{eqnarray*}
&&(\varphi^{(1)}_{t_1}\ast_1 f)(x_1,\,x_2)\\
&&\hs=\int_{\rr^{n_1}}
\lf[\varphi_{t_1}^{(1)}(y_1)-\sum_{|\gz|\le
s_1}\frac{1}{\gz!}(y_1-x_1)^\gz
(D_1^\gz\varphi^{(1)}_{t_1})(x_1)\r]f(x_1-y_1,\,x_2)\,dy_1\\
&&\hs=\int_{|x_1-y_1|<|x_1|/2}\lf[\varphi_{t_1}^{(1)}(y_1)-\sum_{|\gz|\le
s_1}\frac{1}{\gz!}(y_1-x_1)^\gz
(D_1^\gz\varphi^{(1)}_{t_1})(x_1)\r]f(x_1-y_1,\,x_2)\,dy_1\\
&&\hs\hs+\int_{|x_1-y_1|\ge|x_1|/2}\cdots\\
&&\hs\equiv J_1+J_2.
\end{eqnarray*}

On the estimation for $J_1$, notice that if $|z_1|\le
|x_1-y_1|<|x_1|/2$, then $|x_1|/2\le|x_1-z_1|\le 2|x_1|$.  By this
and $\vz^{(1)}\in\cs_0(\rr^{n_1})$, we have
\begin{eqnarray*}
&&\lf|\varphi^{(1)}_{t_1}(y_1)-\sum_{|\gz|\le
s_1}\frac{1}{\gz!}(y_1-x_1)^\gz
(D_1^\gz\varphi^{(1)}_{t_1})(x_1)\r|\\
&&\quad\ls \sup_{|\gz|={s_1+1}}\sup_{|z_1|\le
|x_1-y_1|}|(D_1^\gz\varphi^{(1)}_{t_1})(x_1-z_1)||x_1-y_1|^{s_1+1}\\
&&\quad\ls t_1^{-n_1-s_1-1}\sup_{|z_1|\le
|x_1-y_1|}\lf(1+\frac{|x_1-z_1|}{t_1}\r)^{-\sz_1}|x_1-y_1|^{s_1+1}\\
&&\quad\ls
t_1^{-n_1-s_1-1}\lf(1+\frac{|x_1|}{t_1}\r)^{-\sz_1}|x_1-y_1|^{s_1+1}.
\end{eqnarray*}
Thus, applying
$$|f(x_1-y_1,\,x_2)|\ls (1+|x_1-y_1|)^{-n_1-s_1-2}(1+|x_2|)^{-\sz_2},\leqno(2.5)$$
we further have
\begin{eqnarray*}
|J_1|&&\ls t_1^{-n_1-s_1-1}(1+|x_2|)^{-\sz_2}
\int_{\rr^{n_1}}\lf(1+\frac{|x_1|}{t_1}\r)^{-\sz_1}
\frac{|x_1-y_1|^{s_1+1}}{(1+|x_1-y_1|)^{n_1+s_1+2}}
\,dy_1\\
&&\ls
t_1^{-n_1-s_1-1}\lf(1+\frac{|x_1|}{t_1}\r)^{-\sz_1}(1+|x_2|)^{-\sz_2}.
\end{eqnarray*}
To estimate $J_2$, if $|x_1|>1$ and $\sz_1>0$, using an estimate
similar to (2.5) and the estimation that
$$\lf|\varphi^{(1)}(y_1)-\sum_{|\gz|\le s_1}\frac{1}{\gz!}(y_1-x_1)^\gz
(D_1^\gz\varphi^{(1)}_{t_1})(x_1)\r|\ls
t_1^{-n_1-s_1-1}|x_1-y_1|^{s_1+1},$$ we obtain
\begin{eqnarray*}
|J_2|&&\ls\int_{|y_1-x_1|\ge|x_1|/2} (1+|x_2|)^{-\sz_2}
t_1^{-n_1-s_1-1}
\frac{|x_1-y_1|^{s_1+1}}{(1+|x_1-y_1|)^{\sz_1+n_1+s_1+1}} \,dy_1\\
&&\ls t_1^{-n_1-s_1-1} (1+|x_2|)^{-\sz_2}\int_{|x_1|/2}^\fz
r_1^{-\sz_1-1}
\,dr_1\\
&&\ls t_1^{-n_1-s_2-1}(1+|x_1|)^{-\sz_1}(1+|x_2|)^{-\sz_2},
\end{eqnarray*}
where in the last step, we used the fact that $|x_1|^{-\sz_1}\ls
(1+|x_1|/t_1)^{-\sz_1}$ for $t_1\ge 1$. If $|x_1|\le1$ or
$\sz_1=0$, by (2.5), we then have
\begin{eqnarray*}
|J_2|&&\ls(1+|x_2|)^{-\sz_2}t_1^{-n_1-s_1-1}\int_0^\fz
\frac{r_1^{n_1+s_1}}{(1+r_1)^{n_1+s_1+2}}\,dr_1\\
&&\ls t_1^{-n_1-s_1-1}\lf(1+\frac{|x_1|}{t_1}\r)^{-\sz_1}
(1+|x_2|)^{-\sz_2}.
\end{eqnarray*}
This gives (ii).

To prove (iii), by an argument similar to (i), we obtain that for
all $t_2\in(0,1]$,
$$|(\varphi^{(2)}_{t_2}\ast_2f)(x_1,\,x_2)|\ls
t_2^{u_2+1}(1+|x_1|)^{-\sz_1}(1+|x_2|)^{-\sz_2},\leqno(2.6)$$ where and in
what follows
$(\varphi^{(2)}_{t_2}\ast_2f)(x_1,\,x_2)\equiv\int_{\rr^{n_2}}
\varphi^{(2)}_{t_2}(y_2) f(x_1,x_2-y_2)\,dy_2.$ Thus, if
$|y_1|<|x_1|/2$, then by the mean value theorem, (2.6) and the fact
that $|x_1-z_1|\sim|x_1|$ for $|z_1|\le|x_1|/2$, we have
$$\begin{array}[t]{lll}
&&\lf|(\varphi^{(2)}_{t_2}\ast_2f)(x_1-y_1,\,x_2)-
\dsum_{|\gz|\le u_1}\frac{1}{\gz!}(y_1-x_1)^\gz
\partial^\gz_1(\varphi^{(2)}_{t_2}\ast_2f)(x_1,\,x_2)\r|\\
&&\quad\le |y_1|^{u_1+1}\dsup_{|\gz|=u_1+1}\dsup_{|z_1|\le|x_1|/2}
|(\varphi^{(2)}_{t_2}\ast_2(D^\gz_1f))(x_1-z_1,\,x_2)|\\
&&\quad\ls t_2^{u_2+1}|y_1|^{u_1+1}(1+|x_1|)^{-\sz_1}(1+|x_2|)^{-\sz_2}.
\end{array}\leqno(2.7)$$
If $|y_1|\ge|x_1|/2$, by the mean value theorem and (2.6), we then
have

$$\begin{array}[t]{lll}
&&\lf|(\varphi^{(2)}_{t_2}\ast_2f)(x_1-y_1,\,x_2)-
\dsum_{|\gz|\le s_1}\frac{1}{\gz!}(y_1-x_1)^\gz
\partial^\gz_1(\varphi^{(2)}_{t_2}\ast_2f)(x_1,\,x_2)\r|\\
&&\quad \ls t_2^{u_2+1}|y_1|^{u_1+1}(1+|x_2|)^{-\sz_2}.
\end{array}\leqno(2.8)$$
Noticing that
$$(\varphi_{t_1,\,t_2}\ast f)(x_1,\,x_2)
=(\varphi^{(1)}_{t_1}\ast_1 (\varphi^{(2)}_{t_2}\ast_2
f))(x_1,\,x_2),\leqno(2.9)$$ replacing (2.3) and (2.4) respectively
by (2.7) and (2.8), and repeating the proof of (i), we obtain
(iii).

For (v), by (2.6), we have
$$|(\varphi^{(2)}_{t_2}\ast_2 f)(x_1-y_1,\,x_2)| \ls
t_2^{u_2+1}(1+|x_2|)^{-\sz_2}(1+|x_1-y_1|)^{-n_1-s_1-2}$$ for all
$t_2\in(0,1]$. Replacing (2.5) by this estimate, using (2.9) and repeating
the proof of (ii) lead to (v). A similar argument to (v) yields
(iv).

To obtain (vi), by an  argument similar to (ii), we obtain
$$|(\varphi^{(2)}_{t_2}\ast_2 f)(x_1-y_1,\,x_2)| \ls
t_2^{-n_1-s_1-1}(1+|x_1-y_1|)^{-n_1-s_1-2}\lf(1+\frac{|x_2|}{t_2}\r)^{-\sz_2}$$
for all $t_2\in[1,\fz)$. Replacing (2.5) by this, using (2.9) and
repeating the proof of (ii) leads to (vi). This finishes the
proof of Lemma 2.2.
\end{pf}

\newtheorem{pft}{\it Proof of Lemma 2.1.}
\renewcommand\thepft{}

\begin{pft}\rm
Let $\ez_1\in(0,\,1)$.
Notice that for all $t_1\in(0,\,\fz)$,
$|y_1|\le t_1$ and $x\in\rr^{n_1}$, we have
$t_1+|x_1|\le 2(t_1+|x_1-y_1|).$
By this and Lemma 2.2 (iii) and (iv) ,
we have that for any $t_1\in(0,\ez_1)$, $t_2\in(0,\,1)$,
 $|y_1|<t_1$, $|y_2|<t_2$ and $(x_1,\,x_2)\in\rnt$,
$$|(\varphi_{t_1,\,t_2}\ast f)(x_1-y_1,\,x_2-y_2)|\ls
t_1t_2(1+|x_1|)^{-\sz_1}(1+|x_2|)^{-\sz_2},\leqno(2.10)$$ and that for any
$t_1\in(0,\ez_1]$, $t_2\in[1,\fz)$, $|y_1|<t_1$, $|y_2|<t_2$ and
$(x_1,\,x_2)\in\rnt$,
$$|(\varphi_{t_1,\,t_2}\ast f)(x_1-y_1,\,x_2-y_2)|\ls
t_1t_2^{\sz_2-n_2-s_2-1}(1+|x_1|)^{-\sz_1}(1+|x_2|)^{-\sz_2}.\leqno(2.11)$$
From this and $\sz_2<n_2+s_2+1$, it follows that
\begin{eqnarray*}
&&\sup_{(x_1,\,x_2)\in\rnt}(1+|x_1|)^{\sz_1}(1+|x_2|)^{\sz_2}\int_0^{\ez_1}\int_0^\fz
\int_\rnt|\psi_{t_1,\,t_2}(y_1,\,y_2)|\\
&&\quad\quad\times |(\varphi_{t_1,\,t_2}\ast f)(x_1-y_1,\,x_2-y_2)|
\,dy_1\,dy_2\frac{dt_1}{t_1}\,\frac{dt_2}{t_2}\\
&&\quad\ls \int_0^{\ez_1}\int_0^\fz\int_\rnt\frac1{1+t^{n_2+s_2+2-\sz_2}}
|\varphi_{t_1,\,t_2}(y_1,\,y_2)|\,dy_1\,dy_2\,dt_1\,dt_2\\
&&\quad\ls\ez_1.
\end{eqnarray*}

Let $L_1>1$. By Lemma 2.1 (v) and (vi),  we have that for any
$t_1\in(L_1,\,\fz)$, $t_2\in(0,\,1)$, $|y_1|<t_1$, $|y_2|<t_2$ and
$(x_1,\,x_2)\in\rnt$,
$$|(\varphi_{t_1,\,t_2}\ast f)(x_1-y_1,\,x_2-y_2)|\ls
t_1^{\sz_1-n_1-s_1-1}t_2(1+|x_1|)^{-\sz_1}(1+|x_2|)^{-\sz_2},\leqno(2.12)$$
and that for any $t_1\in(L_1,\,\fz)$, $t_2\in[1,\,\fz)$,
$|y_1|<t_1$, $|y_2|<t_2$ and $(x_1,\,x_2)\in\rnt$,
$$|(\varphi_{t_1,\,t_2}\ast f)(x_1-y_1,\,x_2-y_2)|\ls
t_1^{\sz_1-n_1-s_1-1}t_2^{\sz_2-n_2-s_2-1}(1+|x_1|)^{-\sz_1}
(1+|x_2|)^{-\sz_2}.\leqno(2.13)$$
From this, (2.12), $\sz_1<n_1+s_1+1$ and
$\sz_2<n_2+s_2+1$, it follows that
\begin{eqnarray*}
&&\sup_{(x_1,\,x_2)\in\rnt}(1+|x_1|)^{\sz_1}(1+|x_2|)^{\sz_2}\int_{L_1}^\fz\int_0^\fz
\int_\rnt|\varphi_{t_1,\,t_2}(y_1,\,y_2)|\\
&&\quad\quad\times |(\varphi_{t_1,\,t_2}\ast f)(x_1-y_1,\,x_2-y_2)|
\,dy_1\,dy_2\frac{dt_1}{t_1}\frac{dt_2}{t_2}\\
&&\quad\ls\int_{L_1}^\fz\int_0^\fz
\int_\rnt
|\varphi_{t_1,\,t_2}(y_1,\,y_2)|\,dy_1\,dy_2\,\frac{dt_1}{t_1^{n_1+s_1+2-\sz_1}}
\,\frac{dt_2}{1+t_2^{n_2+s_2+2-\sz_2}} \\
&&\quad\ls (L_1)^{\sz_1-n_1-s_1-1}.
\end{eqnarray*}
Using the symmetry, we then obtain the desired estimates for the
cases $\ez_2\in(0,\,1)$, $L_2\in(1,\,\fz)$,
$(t_1,\,t_2)\in(0,\,\fz)\times (0,\,\ez_2)$ or
$(t_1,\,t_2)\in(0,\,\fz)\times (L_2,\,\fz)$, which gives the first inequality
of Lemma 2.1.

To prove (2.2), notice that if $|y_1|>2L_1>2$ and $|x_1-y_1|<t_1<L_1$, we have
$|x_1|>|y_1|-|x_1-y_1|>L_1$. Then
by (2.10) through (2.13) with $\sz_i$ replaced by
$\sz_i'\in(\sz_i,\,n_1-s_1-1)$, we have
\begin{eqnarray*}
&&\dsup_{(x_1,\,x_2)\in\rnt}(1+|x_1|)^{\sz_1}(1+|x_2|)^{\sz_2}
\int_0^{L_1}\int_0^\fz
\int_{[B^{(1)}(0,\,2L_1)]^\complement\times\rr^{n_2}}
|(\varphi_{t_1,\,t_2}\ast f)(y_1,\,y_2)|\\
&&\quad\quad\quad\times
 |\psi_{t_1,\,t_2}(x_1-y_1,\,x_2-y_2)|
\,dy_1\,dy_2\frac{dt_1}{t_1}\frac{dt_2}{t_2}\\
&&\quad\quad\ls \dsup_{|x_1|>L_1,\,x_2\in\rr^{n_2}}(1+|x_1|)
^{\sz_1-\sz_1'}(1+|x_2|)^{\sz_2-\sz_2'}\int_0^\fz\int_0^\fz
\int_{\rnt}|\psi_{t_1,\,t_2}(y_1,\,y_2)|\\
&&\quad\quad\quad\times
\frac 1{1+t^{n_1-s_1+2-\sz_1'}}
\frac 1{1+t^{n_2-s_2+2-\sz_2'}}
\,dy_1\,dy_2\,dt_1\,dt_2\\
&&\quad\quad\ls (L_1)^{\sz_1-\sz_1'},
\end{eqnarray*}
which gives (2.2) and hence completes the
proof of Lemma 2.1.
\end{pft}

Let $p\in(0,1]$, $s_i\ge\lfloor n_i(1/p-1)\rfloor$ and
$\varphi\in\cs_{s_i}(\rr^{n_i})$ such that (2.1) holds for
$i=1,\,2$. For $f\in\cs'(\rnt)$ and $(x_1,x_2)\in\rnt$, we define
\begin{eqnarray*}
&&S(f)(x_1,\,x_2)\\
&&\quad\equiv\lf(\int_0^\fz\int_0^\fz
\int_{|y_1-x_1|<t_1}\int_{|y_2-x_2|<t_2}|(\varphi_{t_1t_2}\ast
f)(y_1,\,y_2)|^2\,dy_1\,dy_2
\frac{dt_1}{t_1^{n_1+1}}\frac{dt_2}{t_2^{n_2+1}}\r)^{1/2}.
\end{eqnarray*} It is well-known that
$f\in H^p(\rnt)$ if and only if $f\in\cs'(\rnt)$ and $S(f)\in
L^p(\rnt)$. Moreover,
$$\|f\|_{H^p(\rnt)}\sim \|S(f)\|_{L^p(\rnt)};$$ see \cite{cf1,cf2,cf3,fe1}.
Using this fact, Lemma 2.1 and some ideas from
\cite{cf1,cf2,cf3,fe1}, we obtain the following conclusion.

\begin{lem}\hspace{-0.2cm}{.}\label{l2.3}
Let $p\in(0,1]$, $s_i\ge\lfloor n_i(1/p-1)\rfloor$ and
$\sz_i\in(\max\{n_i+s_i,\,n_i/p\},\,n_i+s_i+1)$ for $i=1,\,2$. Then
for any $f\in \cd_{s_1,\,s_2}(\rnt)$, there exist numbers
$\{\lz_k\}_{k\in\nn}\subset\cc$ and $(p,\,2,\,s_1,\,s_2)$-atoms
$\{a_k\}_{k\in\nn}\subset\cd_{s_1,\,s_2}(\rnt)$ such that
$f=\sum_{k\in\nn}\lz_ka_k$ in
$\cd_{s_1,\,s_2;\,\sz_1,\,\sz_2}(\rnt)$ and
$\lf\{\sum_{k\in\nn}|\lz_k|^p\r\}^{1/p}\le C\|f\|_{H^p(\rnt)}$,
where $C$ is a positive constant independent of $f$.
\end{lem}

\begin{pf}\rm We use $\crz$ to denote the set of all dyadic rectangles in $\rnt$.
For $k\in\zz$, let
$$\boz_k\equiv\{(x_1,\,x_2)\in\rnt:\  S(f)(x_1,\,x_2)>2^k\}$$
and
$$\wz \boz_k\equiv\{(x_1,\,x_2)\in\rnt:\
M_s(\chi_{_{\boz_k}})(x_1,\,x_2)>1/2\},$$ where $M_s$ denotes the
strong maximal operator on $\rnt$.
It is easy to see that $\boz_k$ is bounded set. In fact,
observing that $1+|x_i|\le t_i+|x_i|\sim t_i+|y_i|$ for
 $|x_i-y_i|<t_i$ and $t_i\ge1$,
by Lemma 2.2 and $n_i+s_i+1-\sz_i>0$, we have
\begin{eqnarray*}
&&[S(f)(x_1,\,x_2)]^2\\
&&\quad\ls\int_0^1\int_0^1
\int_{|y_1-x_1|<t_1}\int_{|y_2-x_2|<t_2} (1+|y_1|)^{-2\sz_1}
(1+|y_2|)^{-2\sz_2}
\,dy_1\,dy_2\,
\frac{dt_1}{t_1^{n_1}}\,\frac{dt_2}{t_2^{n_2}}\\
&&\quad\quad+\int_1^\fz\int_0^1
\int_{|y_1-x_1|<t_1}\int_{|y_2-x_2|<t_2}
 (1+\frac{|y_1|}{t_1})^{-2\sz_1}(1+|y_2|)^{-2\sz_2}
\,dy_1\,dy_2\,
\frac{dt_1}{t_1^{3n_1+2s_1+3}}\,\frac{dt_2}{t_2^{n_2}}\\
&&\quad\quad+\int_0^1\int_1^\fz
\int_{|y_1-x_1|<t_1}\int_{|y_2-x_2|<t_2} (1+|y_1|)^{-2\sz_1}
(1+\frac{|y_2|}{t_2})^{-2\sz_2}\,dy_1\,dy_2\,
\frac{dt_1}{t_1^{n_1}}\,\frac{dt_2}{t_2^{3n_2+2s_2+3}}\\
&&\quad\quad+\int_1^\fz\int_1^\fz
\int_{|y_1-x_1|<t_1}\int_{|y_2-x_2|<t_2}
(1+\frac{|y_1|}{t_1})^{-2\sz_1}(1+\frac{|y_2|}{t_2})^{-2\sz_2}
\,dy_1\,dy_2\,\\
&&\quad\quad\quad\times
\frac{dt_1}{t_1^{2n_1+s_1+2}}\,\frac{dt_2}{t_2^{3n_2+2s_2+3}}\\
&&\quad\ls (1+|x_1|)^{-2\sz_1}(1+|x_2|)^{-2\sz_2}.
\end{eqnarray*}
Thus for any $k\in\zz$, $\boz_k$ is a bounded set in $\rnt$
and so is $\wz\boz_k$.

For each dyadic rectangle
$R=I\times J$, set
$$\ca(R)\equiv\{(y_1,y_2,t_1,\,t_2):\ (y_1,\,y_2)\in R,\
\sqrt{n_1}|I|<t_1\le 2\sqrt{n_1}|I|,\
\sqrt{n_2}|J|<t_2\le2\sqrt{n_2}|J|\},$$ and
$$\crz_k\equiv\lf\{R\in\crz:\  |R\cap\boz_k|\ge1/2,\ |R\cap
\boz_{k+1}|<1/2\r\}.$$ Obviously, for each $R\in\crz$, there
exists a unique $k\in\zz$ such that $R\in\crz_k$.

From (2.1), for any $(x_1,\,x_2)\in\rnt$, it is easy to see that
$$f(x_1,\,x_2)=\sum_{k\in\zz}\lf\{\sum_{R\in\crz_k}
\int_{\ca(R)}\psi_{t_1,\,t_2}(x_1-y_1,\,x_2-y_2)
(\varphi_{t_1,\,t_2}\ast f)(y_1,\,y_2)\,dy_1\,dy_2
\frac{dt_1}{t_1}\frac{dt_2}{t_2}\r\}.$$ Let
$\lz_k\equiv\frac1C2^k|\boz_k|^{1/p}$ and
$$a_k(x_1,\,x_2)\equiv\lz_k^{-1}
\sum_{R\in\crz_k}\int_{\ca(R)}\psi_{t_1,\,t_2}(x_1-y_1,\,x_2-y_2)
(\varphi_{t_1,\,t_2}\ast f)(y_1,\,y_2)\,dy_1\,dy_2
\frac{dt_1}{t_1}\frac{dt_2}{t_2},$$
where $C$ is a positive constant.
By the argument used in \cite{cf1,cf2,cf3,fe1}, we see that if
we suitably choose the constant $C$, then $\{a_k\}_{k\in\zz}$ are
$(p,\,2,\,s_1,\,s_2)$-atoms and
$$\lf\{\sum_{k\in\zz}|\lz_k|^p\r\}^{1/p}\ls\|f\|_{H^p(\rnt)}.$$

It
remains to prove that $f=\sum_{k\in\zz}\lz_ka_k$ converges in
$\cd_{s_1,\,s_2;\,\sz_1,\,\sz_2}(\rnt)$.
Since $\wz \boz_k$ is bounded, we may assume that
$\wz \boz_k\subset B^{(1)}(0,\, 2^{L_1})\times B^{(2)}(0,\,2^{L_2})$.
Then for any $\az\in\zz_+^{n_1}$
and $\bz\in\zz_+^{n_2}$, by Lemma 2.2, we have
\begin{eqnarray*}
&&\sum_{R\in\crz_k}
\int_{\ca(R)}|(\partial_{x_1}^\az\partial_{x_2}^\bz
\psi_{t_1,\,t_2})(x_1-y_1,\,x_2-y_2)|
(\varphi_{t_1,\,t_2}\ast f)(y_1,\,y_2)|\,dy_1\,dy_2
\frac{dt_1}{t_1}\frac{dt_2}{t_2}\\
&&\quad\ls \sum_{R\in\crz_k} \int_{\ca(R)}
|(\varphi_{t_1,\,t_2}\ast f)(y_1,\,y_2)|\,dy_1\,dy_2
\frac{dt_1}{t_1^{1+|\az|+n_1}}\frac{dt_2}{t_2^{1+|\bz|+n_2}}\\
&&\quad\ls\int_{B^{(1)}(0,\,2^{L_2})}
\int_{B^{(2)}(0,\,2^{L_2})}
\int_0^{L_1}\int_0^{L_2}\,dt_1\,dt_2\,dy_1\,dy_2<\fz,
\end{eqnarray*}
where $(x_1,\,x_2)\in\wz \boz_k$. This shows that
$a_k\in\cd_{s_1,\,s_2,\,\sz_1,\,\sz_2}(\rnt)$. Moreover,
assume that $\supp f\subset
B^{(1)}(0,r_1)\times B^{(2)}(0,r_2)$.
For any $N_i>1+\log r_i$ with $i=1,\,2$, let
$$E_{N_1,\,N_2}\equiv B^{(1)}(0,\,2^{N_1})\times B^{(2)}(0,\,2^{N_2})
\times[2^{-N_1},\,2^{N_1}]\times[2^{-N_2},\,2^{N_2}].$$
Then there exist finite dyadic rectangles $R$, whose set is
denoted by $\crz^{N_1,\,N_2}$, such that
$\ca(R)\cap E_{N_1,\,N_2}\ne\emptyset.$
For each $R\in\crz^{N_1,\,N_2}$, there exists a unique $k\in\zz$ such
that $R\in\boz_k$. Let $K_{N_1,\,N_2}$ be the maximal integer of the absolute
values of all such $k$. Then for $K>K_{N_1,\,N_2}$, by the facts
$\crz^{N_1,\,N_2}\subset \cup_{|k|\le K}\crz_k$ and
Lemma 2.1 together with $\sz_i<\sz_i'<n_i+s_i+1$ for $i=1,2$, we then
have
\begin{eqnarray*}
&&\lf\|f-\sum_{|k|\le K}\lz_ka_k\r\|
_{\cd_{s_1,\,s_2;\,\sz_1,\,\sz_2}(\rnt)}\\
&&\hs\ls\sup_{(x_1,\,x_2)\in\rnt}(1+|x_1|)^{\sz_1}(1+|x_2|)^{\sz_2}\\
&&\quad\quad\times\lf(\int_0^{2^{-N_1}}\int_0^\fz
+\int_{2^{N_1}}^\fz\int_0^\fz+\int_0^\fz\int_0^{
2^{-N_2}}
+\int_0^\fz\int_{2^{N_2}}^\fz\r)\int_\rnt|(\varphi_{t_1,\,t_2}\ast
f)(y_1,\,y_2)|\\
&&\quad\quad\times |\varphi_{t_1,\,t_2}(x_1-y_1,\,x_2-y_2)|
\,dy_1\,dy_2\frac{dt_1}{t_1}\frac{dt_2}{t_2}
+\sup_{(x_1,\,x_2)\in\rnt}(1+|x_1|)^{\sz_1}(1+|x_2|)^{\sz_2}\\
&&\quad\quad\times\int_{2^{-N_1}}^{2^{N_1}}\int_0^\fz
\int_{[B^{(1)}(0,\,2^{N_1})]^\complement \times\rr^{n_2}}|(\varphi_{t_1,\,t_2}\ast
f)(y_1,\,y_2)|\\
&&\quad\quad\times |\psi_{t_1,\,t_2}(x_1-y_1,\,x_2-y_2)|
\,dy_1\,dy_2\frac{dt_1}{t_1}\frac{dt_2}{t_2}
+\sup_{(x_1,\,x_2)\in\rnt}(1+|x_1|)^{\sz_1}(1+|x_2|)^{\sz_2}\\
&&\quad\quad\times
\int_0^\fz\int_{2^{-N_2}}^{2^{N_2}}
\int_{\rr^{n_1}\times[B^{(2)}(0,\,2^{N_2})]^\complement}|(\varphi_{t_1,\,t_2}\ast
f)(y_1,\,y_2)|\\
&&\quad\quad\times |\varphi_{t_1,\,t_2}(x_1-y_1,\,x_2-y_2)|
\,dy_1\,dy_2\frac{dt_1}{t_1}\frac{dt_2}{t_2}\\
&&\quad\ls 2^{-N_1}+2^{-N_2}+2^{N_1(\sz_1-\sz'_1)}+
2^{N_2(\sz_2-\sz'_2)}.
\end{eqnarray*}
This implies the desired conclusion and hence, finishes
the proof of Lemma 2.3.
\end{pf}

The following result plays a key role in the proof of Theorem 1.2.
In what follows, for any $f\in\cd(\rnt)$, we set
$$\sup_{x_2\in{\rr^{n_2}}}\diam(\supp f(\cdot,\,x_2))\equiv
\sup_{x_1,\,y_1\in{\rr^{n_1}},\,x_2\in{\rr^{n_2}}}
\lf\{|x_1-y_1|:\ f(x_1,\,x_2)\ne 0,\ f(y_1,\,x_2)\ne 0\r\},$$ and
$\sup_{x_1\in{\rr^{n_1}}}\diam(\supp f(x_1,\cdot))$ is similarly
defined by interchanging $x_1$ and $x_2$, and $y_1$ and $y_2$.

\begin{lem}\hspace{-0.2cm}{.}\label{2.4}
Let $p\in(0,\,1]$,  $q\in[p,\,1]$ and $\cbq$ be a $q$-quasi-Banach
space. Let $s_1,\,s_2\in\zz_+$ and $T$ be a $\cb_q$-sublinear
operator from $\cd_{s_1,\,s_2}(\rnt)$ to $\cbq$. If there exists a
positive constant $C$ such that for any $f\in\cd_{s_1,\,s_2}(\rnt)$,
\begin{eqnarray*}
\|Tf\|_\cbq&\le& C\lf[\sup_{x_2\in{\rr^{n_2}}}
\diam(\supp f(\cdot,x_2))\r]^{n_1/p}\\
&&\times\lf[\sup_{x_1\in{\rr^{n_1}}} \diam(\supp
f(x_1,\cdot))\r]^{n_2/p}\|f\|_{L^\fz(\rnt)},
\end{eqnarray*}
then $T$ can be extended as a bounded $\cb_q$-sublinear operator from
$\cd_{s_1,\,s_2;\,\sz_1,\,\sz_2}(\rnt)$ to $\cbq$.
\end{lem}

\begin{pf}\rm
Let $\psi\in C^\fz(\rr)$
such that $0\le \psi(x)\le1$ for all $x\in\rr$,
$\psi(x)=1$ if $|x|\le1/2$ and $\psi(x)=0$ if $|x|\ge1$.
Let $\phi(x)\equiv\psi(x/2)-\psi(x)$ for all $x\in\rr$. Then
$\supp\phi\subset \{x\in\rr:\  1/2\le|x|\le2\}$ and
$\sum_{j\in\zz}\phi(2^{-j}x)=1$ for all $x\in\rr\setminus\{0\}$.
Let $\Phi_j(x)\equiv\phi(2^{-j}x)$ for all $x\in\rr$ and $j\in\nn$,
and $\Phi_0(x)\equiv 1-\sum_{j=1}^\fz\phi(2^{-j}x)$ for all $x\in\rr$.
Then $\sum_{j\in\zz_+}\Phi_j(x)=1$ for all $x\in\rr$.

Let $i=1,\, 2$.  For  $j_i\in\zz_+$
and $x_i\in\rr^{n_i}$, let $\Phi^{(i)}_{j_i}(x_i)\equiv\Phi_{j_i}(|x_i|)$.
Then for all
$x_i\in\rr^{n_i}$, we have $\sum_{j_i\in\zz_+}\Phi^{(i)}_{j_i}(x_i)=1$.
Set $R^{(i)}_0\equiv B^{(i)}(0,\,2)$
and $R^{(i)}_{j_i}\equiv\{x_i\in\rr^{n_i}:\ 2^{j_i-1}\le|x_i|\le
2^{j_i+1}\}$ for $j_i\in\nn$. Then
$\supp\Phi^{(i)}_{j_i}\subset R^{(i)}_{j_i}$ for $j_i\in\zz_+$.
For  $j_i\in\zz_+$,
let $\{\wz \psi^{(i)}_{j_i,\,\az_i}:\ |\az_i|\le s_i\}
\subset C^\fz(\rn)$ be the dual basis of $\{x_i^{\az_i}:\  |\az_i|\le s_i\}$
with respect to weight $\Phi^{(i)}_{j_i}|R^{(i)}_{j_i}|^{-1}$, namely, for all
 $\az_i,\,\bz_i\in\zz_+$ with $|\az_i|\le s_i$ and $|\bz_i|\le s_i$,
$$\frac1{|R^{(i)}_{j_i}|}\int_{\rr^{n_i}}
x_i^{\bz_i}\wz \psi^{(i)}_{j_i,\,\az_i}(x_i)\Phi_{j_i}^{(i)}(x_i)\,dx_i=\dz_{\az_i,\bz_i}.$$
Let $\psi^{(i)}_{j_i,\,\az_i}\equiv
|R^{(i)}_{j_i}|^{-1}\wz \psi^{(i)}_{j_i,\,\az_i}\Phi^{(i)}_{j_i}.$
Then for $j_i\in\nn$ and $x_i\in\rr^{n_i}$, we have
$$\psi^{(i)}_{j_i,\,\az_i}(x_i)=2^{-(j_i-1)(n_i+|\az_i|)}
\psi^{(i)}_{1,\,\az_i}(2^{-(j_i-1)}x_i).$$
From this, it is easy to see that for all $j_i\in\zz_+$ and $|\az_i|\le s$,
$$\|\psi^{(i)}_{j_i,\,\az_i}\|_{L^\fz(\rr^{n_i})}\ls 2^{-j_i(n_i+|\az_i|)}.
\leqno(2.14)$$

For $f\in\cd_{s_1,\,s_2}(\rnt)$, assume that $\supp f\subset
B^{(1)}(0,\,2^{k_1})\times B^{(2)}(0, \,2^{k_2})$ for some
$k_1,\,k_2\in\nn$ and $\|f\|_{\cd_{s_1,\,s_2;\,\sz_1,\,\sz_2}(\rnt)}=1$
by the $\cb_q$-sublinear property of $T$.
For $j_1,\,j_2\in\zz_+$, we set
$f_{j_1,\,j_2}\equiv f\Phi^{(1)}_{j_1}\Phi^{(2)}_{j_2} $, and for
any $(x_1,\,x_2)\in\rnt$,
$$P^{(1)}_{j_1,\,j_2}(x_1,\,x_2)
\equiv\sum_{|\az_1|\le
s_1}\psi^{(1)}_{j_1,\,\az_1}(x_1)\int_{\rr^{n_1}} f_{j_1,\,j_2}
(y_1,\,x_2)y_1^{\az_1}\,dy_1,$$
$$P^{(2)}_{j_1,\,j_2}(x_1,\,x_2)
\equiv\sum_{|\az_2|\le s_2}\psi^{(2)}_{j_2,\,\az_2}(x_2)
\int_{\rr^{n_2}} f_{j_1,\,j_2} (x_1,\,y_2)y_2^{\az_2}\,dy_2$$ and
$$P_{j_1,\,j_2}(x_1,\,x_2)
\equiv\sum_{|\az_1|\le s_1}\sum_{|\az_2|\le s_2}
\psi^{(1)}_{j_1,\,\az_1}(x_1)\psi^{(2)}_{j_2,\,\az_2}(x_2)
\int_{\rnt} f_{j_1,\,j_2}
(y_1,\,y_2)y_1^{\az_1}y_2^{\az_2}\,dy_1\,dy_2.$$ Then
\begin{eqnarray*}
f&&=\sum_{j_1=0}^{k_1+1}\sum_{j_2=0}^{k_2+1}\lf(f_{j_1,\,j_2}-
P^{(1)}_{j_1,\,j_2}-P^{(2)}_{j_1,\,j_2}+P_{j_1,\,j_2}\r)\\
&&\quad+\sum_{j_1=0}^{k_1+1}\sum_{j_2=0}^{k_2+1}\lf(P^{(1)}_{j_1,\,j_2}-P_{j_1,\,j_2}\r)
+\sum_{j_1=0}^{k_1+1}\sum_{j_2=0}^{k_2+1}\lf(P^{(2)}_{j_1,\,j_2}-P_{j_1,\,j_2}\r)
+\sum_{j_1=0}^{k_1+1}\sum_{j_2=0}^{k_2+1}P_{j_1,\,j_2}.
\end{eqnarray*}

By the definition of $\cd_{s_1,\,s_2;\,\sz_1,\,\sz_2}(\rnt)$, it
is easy to see that
$$\|f_{j_1,\,j_2}\|_{L^\fz(\rnt)}\ls
2^{-j_1\sz_1}2^{-j_2\sz_2}.\leqno(2.15)$$
Using $ \|\Phi^{(i)}_{j_i}
\|_{L^\fz(\rr^{n_i})}\le 1$, we obtain
$$\lf\|\dint_{\rr^{n_1}} f_{j_1,\,j_2}(y_1,\,\cdot)
y_1^{\az_1}\,dy_1\r\|_{L^\fz(\rr^{n_2})} \ls
2^{j_1(n_1+|\az_1|-\sz_1)}2^{-j_2\sz_2}, \leqno(2.16)$$
$$\lf\|\dint_{\rr^{n_2}} f_{j_1,\,j_2}(\cdot,\,y_2)y_2^{\az_2}\,dy_2\r\|_{L^\fz(\rr^{n_1})}
\ls
2^{-j_1\sz_1}2^{j_2(n_2+|\az_2|-\sz_2)},
\leqno(2.17)$$
and
$$\begin{array}[t]{l}
\lf|\dint_{\rnt} f_{j_1,\,j_2}(y_1,\,y_2)
y_1^{\az_1}y_2^{\az_2}\,dy_1\,dy_2\r|
\ls 2^{j_1(n_1+|\az_1|-\sz_1)}2^{j_2(n_2+|\az_2|-\sz_2)}.
\end{array}\leqno(2.18)$$
By the estimates (2.14) through (2.18), we have
$$\lf\|f_{j_1,\,j_2}-
P^{(1)}_{j_1,\,j_2}-P^{(2)}_{j_1,\,j_2}+P_{j_1,\,j_2}\r\|_{L^\fz(\rnt)}\\
\ls
2^{-j_1\sz_1}2^{-j_2\sz_2}.$$
Since $f_{j_1,\,j_2}-
P^{(1)}_{j_1,\,j_2}-P^{(2)}_{j_1,\,j_2}+P_{j_1,\,j_2}\in\cd_{s_1,\,s_2}(\rnt),$
by the assumption of the lemma, we then have
$$\lf\|T\lf(f_{j_1,\,j_2}-
P^{(1)}_{j_1,\,j_2}-P^{(2)}_{j_1,\,j_2}+P_{j_1,\,j_2}\r)\r\|_{\cbq}\ls
2^{j_1(n_1/p-\sz_1)}2^{j_2(n_2/p-\sz_2)},$$
and hence, by $\sz_i>n_i/p$ for $i=1,\,2$,
$$\begin{array}[t]{ccl}
&&\lf\|T\lf[\dsum_{j_1=0}^{k_1+1}\dsum_{j_2=0}^{k_2+1}\lf(f_{j_1,\,j_2}-
P^{(1)}_{j_1,\,j_2}-P^{(2)}_{j_1,\,j_2}+P_{j_1,\,j_2}\r)\r]\r\|_{\cbq}\\
&&\quad\ls \lf\{\dsum_{j_1=0}^{k_1+1}\dsum_{j_2=0}^{k_2+1}
2^{j_1q(n_1/p-\sz_1)}2^{j_2q(n_2/p-\sz_2)}\r\}^{1/q}\ls1.
\end{array}\leqno(2.19)$$

Moreover, we write
\begin{eqnarray*}
&&\sum_{j_1=0}^{k_1+1}\sum_{j_2=0}^{k_2+1}\lf[P^{(1)}_{j_1,\,j_2}(x_1,\,x_2)
-P_{j_1,\,j_2}(x_1,\,x_2)\r]\\
&&\quad=\sum_{|\az_1|\le s_1} \sum_{j_1=1}^{k_1+1}
\sum_{j_2=0}^{k_2+1}\sum_{\ell_1=j_1}^{k_1+1}
\lf[\psi^{(1)}_{j_1,\,\az_1}(x_1)-\psi^{(1)}_{j_1-1,\,\az_1}(x_1)\r]
\lf[\int_{\rr^{n_1}} f_{\ell_1,\,j_2}
(y_1,\,x_2)y_1^{\az_1}\,dy_1\r.\\
&&\quad\quad-\sum_{|\az_2|\le s_2}
\psi^{(2)}_{j_2,\,\az_2}(x_2)\int_{\rr^{n_1}}\int_{{\rr^{n_2}}}
f_{\ell_1,\,j_2} (y_1,\,y_2)y_1^{\az_1}y_2^{\az_2}\,dy_1\,dy_2
\bigg]\\
&&\quad\equiv\sum_{|\az_1|\le s_1} \sum_{j_1=1}^{k_1+1}
\sum_{j_2=0}^{k_2+1}\sum_{\ell_1=j_1}^{k_1+1}A_{\az_1,\,j_1,\,\ell_1,\,j_2}(x_1,\,x_2).
\end{eqnarray*}

By (2.14), (2.15) and (2.18), we have
$$\|A_{\az_1,\,j_1,\,\ell_1,\,j_2}\|_{L^\fz(\rnt)}\ls
2^{-j_1(n_1+|\az_1|)}2^{\ell_1(n_1+|\az_1|-\sz_1)}
2^{-j_2\sz_2}.$$
Noticing that
$A_{\az_1,\,j_1,\,\ell_1,\,j_2}\in\cd_{s_1,\,s_2}(\rnt)$, by the
assumption of the lemma, we obtain
$$\lf\|T\lf(A_{\az_1,\,j_1,\,\ell_1,\,j_2}\r)\r\|_\cbq
\ls2^{j_1(n_1/p-n_1-|\az_1|)}2^{\ell_1(n_1+|\az_1|-\sz_1)}
2^{j_2(n_2/p-\sz_2)}.$$
Thus by $\sz_i\in(\max\{n_i/p,\,n_i+s_i\},\,n_i+s_i+1)$ for
$i=1,\,2$, we further have
$$\begin{array}[t]{ccl}
&&\lf\|T\lf[\dsum_{j_1=0}^{k_1+1}\dsum_{j_2=0}^{k_2+1}
\lf(P^{(1)}_{j_1,\,j_2}-P_{j_1,\,j_2}\r)\r]\r\|_\cbq\\
&&\quad\ls \lf\{\dsum_{|\az_1|\le s_1}
\dsum_{j_1=1}^{k_1+1}
\dsum_{j_2=0}^{k_2+1}\dsum_{\ell_1=j_1}^{k_1+1}
2^{j_1q(n_1/p-n_1-|\az_1|)}2^{\ell_1q(n_1+|\az_1|-\sz_1)}
2^{j_2q(n_2/p-\sz_2)}
\r\}^{1/q}\ls1.
\end{array}\leqno(2.20)$$

Similarly, by symmetry, we have
$$\lf\|T\lf[\sum_{j_1=0}^{k_1+1}\sum_{j_2=0}^{k_2+1}
\lf(P^{(2)}_{j_1,\,j_2}-P_{j_1,\,j_2}\r)\r]\r\|_\cbq
\ls1. \leqno(2.21)$$

Finally, we write
\begin{eqnarray*}
\sum_{j_1=0}^{k_1+1}\sum_{j_2=0}^{k_2+1}P_{j_1,\,j_2}
&&=\sum_{|\az_1|\le s_1}\sum_{|\az_2|\le s_2}
\sum_{j_1=1}^{k_1+1}\sum_{j_2=1}^{k_2+1}
\sum_{\ell_1=j_1}^{k_1+1}\sum_{\ell_2=j_2}^{k_2+1}
\lf(\psi^{(1)}_{j_1,\,\az_1}-\psi^{(1)}_{j_1-1,\,\az_1}\r)\\
&&\quad\quad\times
\lf(\psi^{(2)}_{j_2,\,\az_2}-\psi^{(2)}_{j_2-1,\,\az_2}\r)
\int_{\rr^{n_1}}\int_{\rr^{n_2}}
f_{\ell_1,\,\ell_2}
(y_1,\,y_2)y_1^{\az_1}y_2^{\az_2}\,dy_1\,dy_2\\
&&\quad\equiv\sum_{|\az_1|\le s_1}\sum_{|\az_2|\le s_2}
\sum_{j_1=1}^{k_1+1}\sum_{j_2=1}^{k_2+1}
\sum_{\ell_1=j_1}^{k_1+1}\sum_{\ell_2=j_2}^{k_2+1}
A_{\az_1,\,j_1,\,\ell_1,\,\az_2,\,j_2,\,\ell_2}.
\end{eqnarray*}
From (2.14) and (2.17), it follows that
\begin{eqnarray*}
&&\lf\|A_{\az_1,\,j_1,\,\ell_1,\,\az_2,\,j_2,\,\ell_2}\r\|_{L^\fz(\rnt)}
\ls 2^{-j_1(n_1+|\az_1|)}2^{\ell_1(n_1+|\az_1|-\sz_1)}
2^{-j_2(n_2+|\az_2|)}2^{\ell_2(n_2+|\az_2|-\sz_2)}.
\end{eqnarray*}
Since
$A_{\az_1,\,j_1,\,\ell_1,\,\az_2,\,j_2,\,\ell_2}\in\cd_{s_1,\,s_2}(\rnt)$,
by the assumption of the lemma, then
\begin{eqnarray*}
&&\lf\|T(A_{\az_1,\,j_1,\,\ell_1,\,\az_2,\,j_2,\,\ell_2})\r\|_\cbq\\
&&\quad\ls2^{j_1(n_1/p-n_1-|\az_1|)}2^{\ell_1(n_1+|\az_1|-\sz_1)}
2^{j_2(n_2/p-n_2-|\az_2|)}2^{\ell_2(n_2+|\az_2|-\sz_2)}.
\end{eqnarray*}
From this and $\sz_i\in(\max\{n_i/p,\,n_i+s_i\},\,n_i+s_i+1)$ for
$i=1,\,2$, it follows that
\begin{eqnarray*}
&&\lf\|T\lf(\sum_{j_1=0}^{k_1+1}\sum_{j_2=0}^{k_2+1}
P_{j_1,\,j_2}\r)\r\|_\cbq\\
&&\quad\ls\lf\{\sum_{|\az_1|\le s_1} \sum_{|\az_2|\le s_2}
\sum_{j_1=1}^{k_1+1} \sum_{j_2=1}^{k_2+1}
\sum_{\ell_1=j_1}^{k_1+1} \sum_{\ell_2=j_2}^{k_2+1}
2^{j_1q(n_1/p-n_1-|\az_1|)}2^{\ell_1q(n_1+|\az_1|-\sz_1)}\r.\\
&&\quad\quad\times2^{j_2q(n_2/p-n_2-|\az_2|)}2^{\ell_2q(n_2+|\az_2|-\sz_2)}
\Bigg\}^{1/q}\ls1.
\end{eqnarray*}
By this together with the estimates (2.19) through (2.21) and the
$\cb_q$-sublinear property of $T$, we obtain that
$\|Tf\|_\cbq\ls\|f\|_{\cd_{s_1,\,s_2;\,\sz_1,\,\sz_2}(\rnt)}$,
which implies that $T$ is bounded from $\cd_{s_1,\,s_2;\,
\sz_1,\,\sz_2}(\rnt)$ to $\cbq$. This finishes the proof of Lemma
2.4.
\end{pf}

\newtheorem{pftt}{\it Proof of Theorem 1.1.}
\renewcommand\thepftt{}

\begin{pftt}\rm
The necessity is obvious. In fact, if $T$ extends to a bounded
$\cb_q$-sublinear operator from $H^p(\rnt)$ to $\cbq$, then for any
$(p,\,2,\,s_1,\,s_2)$-atom $a$,
$$\|Ta\|_\cbq\ls \|a\|_{H^p(\rnt)}\ls1.$$

To prove the sufficiency, for any $f\in \cd_{s_1,\,s_2}(\rnt)$,
let
$$\ell_1\equiv\sup_{x_2\in{\rr^{n_2}}}\diam(\supp f(\cdot,\,x_2))$$
and $\ell_2\equiv\sup_{x_1\in{\rr^{n_1}}}\diam(\supp
f(x_1,\,\cdot))$. Then there exists a positive constant $C$
independent of $f$ such that
$C(\ell_1)^{-n_1/p}(\ell_2)^{-n_2/p}\|f\|_{L^\fz(\rnt)}^{-1}f$ is a
$(p,\,2,\,s_1,\,s_2)$-atom, and thus, by the assumption of the
theorem,
$$\|Tf\|_\cbq\ls(\ell_1)^{n_1/p}(\ell_2)^{n_2/p}\|f\|_{L^\fz(\rnt)},$$
which shows that $T$ satisfies the assumptions of Lemma 2.4. For
$i=1,\,2$, choose $\sz_i\in
(\max\{n_i+s_i,\,n_i/p\},\,n_i+s_i+1)$. By Lemma 2.4, $T$ is
bounded from $\cd_{s_1,\,s_2;\,\sz_1,\,\sz_2}(\rnt)$ to $\cbq$.

On the other hand, for any $f\in\cd_{s_1,\,s_2}(\rnt)$, by Lemma
2.3, there exist numbers $\{\lz_j\}_{j\in\nn}\subset\cc$ and
$(p,\,2,\,s_1,\,s_2)$-atoms
$\{a_j\}_{j\in\nn}\subset\cd_{s_1,\,s_2}(\rnt)$ such that
$f=\sum_{j\in\nn}\lz_ja_j$ in
$\cd_{s_1,\,s_2;\,\sz_1,\,\sz_2}(\rnt)$ and
$\{\sum_{j\in\nn}|\lz_j|^p\}^{1/p}\ls \|f\|_{H^p(\rnt)}$. From
this and Lemma 2.4, it follows that $Tf=\sum_{j\in\nn}\lz_jTa_j$ in
$\cbq$. Thus $Tf\in\cbq$ and by the monotonicity of the
sequence space $\ell^q$,
\begin{eqnarray*}
\|Tf\|_\cbq&&\le\lf\{\sum_{j\in\nn}|\lz_j|^q
\|Ta_j\|^q_\cbq\r\}^{1/q} \ls\lf\{\sum_{j\in\nn}|\lz_j|^p
\r\}^{1/p} \ls\|f\|_{H^p(\rnt)}.
\end{eqnarray*}
This together with the density of $\cd_{s_1,\,s_2}(\rnt)$ in
$H^p(\rnt)$ implies that $T$ can be extended as a bounded $\cb_q$-sublinear
operator from $H^p(\rnt)$ to $\cbq$, which completes the proof of
Theorem 1.1.
\end{pftt}

Using Theorem 1.1, we can now prove Corollary 1.1.

\newtheorem{pfc}{\it Proof of Corollary 1.1.}
\renewcommand\thepfc{}

\begin{pfc}\rm
By Theorem 1.1, it suffices to prove that for all smooth atoms $a$,
$\|T(a)\|_{L^{q_0}(\rnt)}\ls1$. To prove this, we follow the
procedure used in the proof of Theorem 1 in \cite{fe1} (see also
\cite{fe2}). Assume that $a$ is a smooth $(p,\,2,\,s_1,\,s_2)$-atom
supported in open set $\boz$. Let
$\wz\boz\equiv\{(x_1,\,x_2)\in\rnt:\ M_s(\chi_\boz)(x_1,\,x_2)>1/2\}$
and
$$\boz_0\equiv\{(x_1,\,x_2)\in\rnt:\ M_s(\chi_{\wz{\wz\boz}})(x_1,\,x_2)>1/16\}.$$
Then $|\boz_0|+|\wz \boz|\ls|\boz|$. By the boundedness of $T$ from
$L^2(\rnt)$ to $L^{q_0}(\rnt)$ and the H\"older inequality, we have
\begin{eqnarray*}
\lf\{\int_{\boz_0}|T(a)(x_1,\,x_2)|^{q}\,dx_1\,dx_2\r\}^{1/q}
&&\ls\lf\{\int_{\boz_0}|T(a)(x_1,\,x_2)|^{q_0}\,dx_1\,dx_2\r\}^{1/q_0}
|\boz|^{1/q-1/q_0}\\
&&\ls\|a\|_{L^2(\rnt)} |\boz|^{1/p-1/2}\ls1.
\end{eqnarray*}
We still need to prove that
$\int_{(\boz_0)^\complement}|T(a)(x_1,\,x_2)|^{q}\,dx_1\,dx_2\ls1.$
Without loss of generality,  we may assume that $q\le1$. The proof
of the case $q\in (1,2)$ is similar and we omit the details. To this end, for each
$R\in \cm(\boz)$, assume that $R=I\times J$. Denote by
$\cm^{(1)}(\wz\boz)$ the set of all maximal subrectangles in the first
direction in $\boz$. Let $\hat R\equiv\hat I\times
J\in\cm^{(1)}(\wz\boz)$ and $\hat{\hat R}\equiv\hat I\times \hat
J\in\cm^{(1)}(\wz{\wz \boz})$, and define
$\gz_1(R,\,\boz)\equiv|\hat I|/|I|$ and $\gz_2(\hat
R,\,\wz\boz)\equiv|\hat J|/|J|$. Then $16\hat{\hat R}\subset\boz_0$.
Notice that by the Journ\'e covering lemma (see \cite{p}), for any
fixed $\dz'>0$,  we have
$$\sum_{R\in \cm(\boz)}[\gz_1(R,\,\boz)]^{-\dz'}
|R|\ls |\boz|\leqno{(2.22)}$$
and
$$\sum_{\hat R\in \cm^{(1)}(\wz\boz)}[\gz_2(\hat R,\,\wz\boz)]^{-\dz'}
|R|\ls |\boz|.\leqno{(2.23)}$$

Since $q\le 1$, we write
\begin{eqnarray*}
&&\int_{(\boz_0)^\complement}|T(a_R)(x_1,\,x_2)|^{q}\,dx_1\,dx_2\\
&&\quad\le
\sum_{R\in \cm(\boz)} \int_{(\boz_0)^\complement}|T(a_R)(x_1,\,x_2)|^{q}\,dx_1\,dx_2\\
&&\quad\le\sum_{R\in \cm(\boz)}\int_{(\rr^{n_1}\setminus 16\hat I)\times
\rr^{n_2}} |T(a_R)(x_1,\,x_2)|^q\,dx_1\,dx_2+ \sum_{R\in \cm(\boz)}
\int_{\rr^{n_1}\times(\rr^{n_2}\setminus16\hat J)} \cdots\\
&&\quad\equiv L_1+L_2.
 \end{eqnarray*}
Noticing that  $a_R|R|^{1/2-1/p}\|a_R\|_{L^2(\rnt)}^{-1}$ is a rectangle atom,
we have
\begin{eqnarray*}
\int_{(\rr^{n_1}\setminus 16\hat I)\times \rr^{m_2}}
|T(a_R)(x_1,\,x_2)|^q\,\,dx_1\,dx_2\ls [\gz_1(R,\,\boz)]^{-\dz}
|R|^{1-q/q_0}\|a_R\|^q_{L^2(\rnt)}.
\end{eqnarray*}
By $1/q_0-1/q=1/2-1/p$ and $p\le1$ and (2.22), we obtain
\begin{eqnarray*}
L_1
&&\ls \lf\{\sum_{R\in \cm(\boz)}\|a_R\|^2_{L^2(\rnt)}\r\}^{q/2}\\
&&\quad\times\lf\{\sum_{R\in \cm(\boz)}[\gz_1(R,\,\boz)]^{-2\dz/(2-q)}
|R|^{[2(q_0-q)]/[q_0(2-q)]}\r\}^{1-q/2} \\
&&\ls|\boz|^{q(1/2-1/p)}|\boz|^{q(1/2-1/q_0)}
\lf\{\sum_{R\in \cm(\boz)}[\gz_1(R,\,\boz)]^{-2\dz/(2-q)}
|R|\r\}^{1-q/2}\\
&&\ls|\boz|^{q(1/2-1/q)}|\boz|^{1-q/2}\ls1.
\end{eqnarray*}
Similarly, by (2.23), we have $L_2\ls1$.
This finishes the proof of Corollary 1.1.
\end{pfc}

\section{\hspace{-0.55cm}{\bf .} Proofs of Theorem 1.2 and Theorem 1.3\label{s3}}

\hskip\parindent To prove Theorem 1.2,
we recall the well-known boundedness of fractional integrals on
$\rn$; see \cite[p.\,117]{st}.

\begin{lem}\hspace{-0.2cm}{.}\label{l3.1}
Let $\az\in (0, 1)$,\ $p\in (1,n/\az)$ and $1/q=1/p-\az/n.$
Let $I_\az$ be the fractional integral operator on $\rn$ defined by
$$I_\az (f)(x)=\dint_\rn \dfrac{f(y)}{|x-y|^{n-\az}}\,dy$$
for $f\in L_\loc^1(\rn)$ and $x\in\rn$. Then $I_\az$ is bounded from
$L^p(\rn)$ to $L^q(\rn)$, namely, there exists a positive constant
$C$ such that for all $f\in L^p(\rn)$,
$$\|I_\az (f)\|_{L^q(\rn)}\le C\|f\|_{L^p(\rn)}.$$
\end{lem}

\newtheorem{pfo}{\it Proof of Theorem 1.2.}
\renewcommand\thepfo{}

\begin{pfo}\rm
Since $[b,\,T]$ is linear with respect to $b$ and $T$, then it
suffices to prove Theorem 1.2 for $b\in\laa$ with $\|b\|_\laa=1$ and
$T$ with $\|K\|=1$. By (K1) and Definition 1.4, we have
$$\begin{array}{cl}
|[b,\,T] (f)(x_1,\,x_2)|
&\ls\dint_{\rnm}\dfrac 1{|x_1- y_1|^{n-\az_1}}
\dfrac 1{|x_2- y_2|^{m-\az_2}}|f(y_1,\,y_2)|\,dy_1\, dy_2\\
&\ls I^{(1)}_{\az_1}\lf[I^{(2)}_{\az_2}(|f|)\r](x_1,\,x_2),
\end{array}$$
where $ I^{(1)}_{\az_1}$ and $I^{(2)}_{\az_2}$ are the fractional
integral operators with respect to $x_1$ or $x_2$, respectively.
By Lemma
3.1, for all $f\in\lp$, we have
$$\begin{array}{cl} \lf\|[b,\,T](f)\r\|_{\lq}
 &\ls\lf\|\lf\|I^{(1)}_{\az_1}\lf[I^{(2)}_{\az_2}(|f|)\r]
 \r\|_{L^q(\rr^m,\,dx_2)}\r\|_{L^q(\rn,\,dx_1)}\\
&\ls\lf\|\lf\|I^{(2)}_{\az_2}(|f|)\r\|_{L^q(\rr^m,\,dx_2)}\r\|_{L^p(\rn,\,dx_1)}\\
&\ls\|f\|_\lp,
\end{array}$$
where and in the sequel, we use $\|\cdot\|_{L^p(\rn,\, dx_1)}$ and
$\|\cdot\|_{L^p(\rr^m,\, dx_2)}$ to
denote the $L^p(\rn)$-norm with respect to the variable $x_1$ and $x_2$
respectively. This
finishes the proof of Theorem 1.2.
\end{pfo}

\newtheorem{pfs}{\it Proof of Theorem 1.3.}
\renewcommand\thepfs{}

\begin{pfs}\rm
Since $[b,\,T]$ is linear with respect to $b$ and $T$, then it suffices to
prove Theorem 1.3 for $b\in\laa$ with $\|b\|_\laa=1$ and $T$ with $\|K\|=1$.
By Theorem 1.1 and Corollary 1.1,
it suffices to prove that there exists a positive $\dz$ such that for all
rectangular $(p,\,2,\,s_1,\,s_2)$-atoms $a$ supported on $R=I\times J$ and
$\gz\ge8\max\{n^{1/2},m^{1/2}\}$,
$$\int_{(\rnm)\setminus \wz R_\gz} |[b,\,T](a)(x_1,\,x_2)|^q\,dx_1\,dx_2
\ls\gz^{-\dz}.\leqno(3.1)$$
Without loss of generality, we may assume
that $R=I\times J=[0,1]^n\times[0,1]^m$. In fact, if
letting $b_{x_1^0,\,x_2^0,\,\ell_1,\ell_2}(x_1,\,x_2)
\equiv\ell_1^{-\az_1}\ell_2^{-\az_2}b(x_1^0+\ell_1x_1,\,
x_2^0+\ell_2x_2),$
$$K_{x_1^0,\,x_2^0,\,\ell_1,\ell_2}(x_1,\,y_1,\,x_2,\,y_2)=
\ell_1^n\ell_2^mK(x_1^0+\ell_1x_1,\, x_1^0+\ell_1y_1,
x_2^0+\ell_2x_2,\,x_2^0+\ell_1y_2)$$
and $T_{x_1^0,\,x_2^0,\,\ell_1,\,\ell_2}$ be a Calder\'on-Zygmund
operator with kernel $K_{x_1^0,\,x_2^0,\,\ell_1,\,\ell_2}$ for
some $x_1^0\in\rn,x_2^0\in\rr^m$ and some $\ell_1,\,\ell_2>0$,
then it is easy to check that
$$\|b_{x_1^0,\,x_2^0,\,\ell_1,\,\ell_2}\|_\laa=\|b\|_\laa=1$$
and $K_{x_1^0,\,x_2^0,\,\ell_1,\,\ell_2}$ also satisfies (K1) through (K4)
with $\|K_{x_1^0,\,x_2^0,\,\ell_1,\,\ell_2}\|=\|K\|=1$.
Moreover, if let $\wz a$ be a rectangular $(p,\,2,\,s_1,\,s_2)$-atom
supported in $R'=I'\times J'=
\{x_1^0+\ell_1 I\}\times \{x_2^0+\ell_2 J\}$, and
 $a(x_1,\,x_2)\equiv\ell_1^n\ell_2^m \wz a(x_1^0+\ell_1x_1,\,x_2^0+\ell_2x_2)$,
then $a$ is a rectangular $(p,\,2,\,s_1,\,s_2)$-atom supported in
$R=[0,\,1]^n\times[0,\,1]^m$, where $x_0^1+\ell_1I=
\{x_0^1+\ell_1x_1: x_1\in I\}$ and
$x_2^0+\ell_2J$ is similarly defined.
By setting $x_i'=x_i^0+\ell_ix_i$ and
 $y_i'=y_i^0+\ell_iy_i$ for $i=1,\,2$, we have
\begin{eqnarray*}
&&[b,\, T](\wz a)(x_1',x_2')\\
&&\quad=\dint_{\rnm} K(x_1',\,y_1',\,x_2',\,y_2')\\
&&\quad\quad\times[b(x_1',\,x_2')
-b(x_1',\,x_2')-b(y_1',\,y_2')+b(y_1',\,y_2')]a'(y_1',\,y_2')\, dy_1'\,dy_2'\\
&&\quad=\ell_1^{\az_1-n}\ell_2^{\az_2-m}
\dint_{\rnm} K_{x_1^0,\,x_2^0,\,\ell_1,\,\ell_2}(x_1,\,y_1,\,x_2,\,y_2)
[b_{x_1^0,\,x_2^0,\,\ell_1,\,\ell_2}(x_1,\,x_2)\\
&&\quad\quad
-b_{x_1^0,\,x_2^0,\,\ell_1,\,\ell_2}(x_1,\,y_2)-
b_{x_1^0,\,x_2^0,\,\ell_1,\,\ell_2}(y_1,\,x_2)+
b_{x_1^0,\,x_2^0,\,\ell_1,\,\ell_2}(y_1,\,y_2)]a(y_1,\,y_2)\, dy_1\,dy_2\\
&&\quad=\ell_1^{\az_1-n}\ell_2^{\az_2-m}
[b_{x_1^0,\,x_2^0,\,\ell_1,\,\ell_2},\,
T_{x_1^0,\,x_2^0,\,\ell_1,\,\ell_2}](a)(x_1,\,x_2),
\end{eqnarray*}
which together with $1/q=1-\az_1/n=1-\az_2/m$ yields
\begin{eqnarray*}
&&\int_{(\rnm)\setminus\wz R'_\gz}
|[b,\,T](\wz a)(x_1',\,x_2')|^q\,dx_1'\,dx_2'\\
&&\quad=\ell_1^n\ell_2^m
\int_{(\rnm)\setminus \wz R_\gz}
|[b,\,T](\wz a)(x_1',\,x_2')|^q\,dx_1'\,dx_2'\\
&&\quad=
\int_{(\rnm)\setminus \wz R_\gz}
|[b_{x_1^0,\,x_2^0,\,\ell_1,\,\ell_2},\,T_{x_1^0,\,x_2^0,\,\ell_1,\,\ell_2}](a)
(x_1,\,x_2)|^q\,dx_1\,dx_2,
\end{eqnarray*}
where $\wz R'$ denotes the $\gz$ fold enlargement of $R'$. Then by
this, (1.3) and the facts that $K_{x_1^0,\,x_2^0,\,\ell_1,\,\ell_2}$
and $b_{x_1^0,\,x_2^0,\,\ell_1,\,\ell_2}$ satisfy the same
conditions as $K$ and $b$ respectively, we may assume that
$R=I\times J=[0,\,1]^n\times[0,\,1]^m$.

Let $a$ be a rectangular $(p,\,2,\,s_1,\,s_2)$-atom supported in
$R=I\times J=[0,\,1]^n\times[0,\,1]^m$. Let $\gz_1\equiv 8n^{1/2}$,
$\gz_2\equiv 8m^{1/2}$ and $\gz\ge\max\{\gz_1,\,\gz_2\}$. Then
\begin{eqnarray*}
&&\int_{(\rnm)\setminus \wz R_\gz} |[b,\,T](a)(x_1,\,x_2)|^q\,dx_1dx_2\\
&&\quad\le\dint_{x_1\not\in \gz I}\dint_{x_2\in \gz_2
J}|[b,\,T](a)(x_1,\,x_2)|^q\,dx_1dx_2+\dint_{x_1\not\in \gz
I}\dint_{x_2\not\in \gz_2J}\cdots+\dint_{x_1\in \gz_1I}
\dint_{x_2\not\in \gz
J}\cdots\\
&&\quad\equiv G_1+G_2+G_3.
\end{eqnarray*}
By symmetry, it suffices to estimate $G_1$ and $G_2$.

The H\"older inequality implies that
$$G_1\ls\dint_{x_1\not\in \gz I}\lf\|[b,\,T]
a(x_1,\,\cdot)\r\|^q_{L^{q_1}(\rr^m,\,dx_2)}dx_1.$$
By $\int_\rn a(x_1,\,x_2)\,dx_1=0$ for all $x_2\in\rr^m$, we have
\begin{eqnarray*}
&&[b,\,T](a)(x_1,\,x_2)\\
&&\quad=\dint_{\rnm} [K(x_1,\,y_1,\,x_2,\,y_2)-K(x_1,\,0,\,x_2,\,y_2)]\\
&&\quad\quad\times[b(x_1,\,x_2)-b(x_1,\,x_2)
-b(x_1,\,y_2)+b(y_1,\,y_2)]a(y_1,\,y_2)\,dy_1\, dy_2\\
&&\quad\quad+\dint_{\rnm} K(x_1,\,0,\,x_2,\,y_2)\\
&&\quad\quad\times[b(0,\,x_2)-b(0,\,y_2)
-b(y_1,\,x_2)+b(y_1,\,y_2)]a(y_1,\,y_2)\,dy_1\, dy_2\\
&&\quad\equiv L_1+L_2.
\end{eqnarray*}
Notice that if $x_1\not\in \gz I$ and  $y_1\in I$,  then $|y_1|\le
|x_1|/2$ and $|x_1-y_1|\ls 2|x_1|$. Thus for any  $x_1\not\in \gz
I$ and $x_2\in\rr^m$,
 by Definition 1.4, (K1), (K2) and the H\"older inequality, we obtain
\begin{eqnarray*}
|L_1|&&\ls\dint_{I}\dint_{J}
\dfrac{|y_1|^{\ez_1}}{|x_1|^{n+\ez_1-\az_1}}\dfrac{1}{|x_2-y_2|^{m-\az_2}}
|a(y_1,\,y_2)|dy_1\, dy_2\\
&&\ls\dfrac{1}{|x_1|^{n+\ez_1-\az_1}}
\dint_{J}\dfrac{1}{|x_2- y_2|^{m-\az_2}}
\lf(\dint_{I}|a(y_1,\,y_2)|^2 d y_1\r)^{1/2}d y_2\\
&&\ls\dfrac{1}{|x_1|^{n+\ez_1-\az_1}}I^{(2)}_{\az_2}
\lf[\|a\|_{L^2(\rn,\,d y_1)}\r](x_2)
\end{eqnarray*}
and
\begin{eqnarray*} |L_2|&&\ls
\dint_{I}\dint_{J}\dfrac{|y_1|^{\az_1}}{|x_1|^n}\dfrac{1}{|x_2-
y_2|^{m-\az_2}}
|a(y_1,\,y_2)|dy_1\, dy_2\\
&&\ls\dfrac{1}{|x_1|^n}I^{(2)}_{\az_2}\lf[
\|a\|_{L^2(\rr^n,\,d y_1)}\r](x_2).
\end{eqnarray*}
Since (1.4) implies that $n-(n+\ez_1-\az_1)q<0$ and $n-nq<0$, then
by (R3) and Lemma 3.1, we obtain
\begin{eqnarray*}
G_1&&\ls\dint_{x_1\not\in \gz
I}\lf(\|L_1\|^q_{L^{q_1}(\rr^m,\,dx_2)}+
\|L_2\|^q_{L^{q_1}(\rr^m,\,dx_2)}\r)
dx_1\\
&&\ls\dint_{x_1\not\in
\gz I}\lf(\dfrac{1}{|x_1|^{(n+\ez_1-\az_1)q}}+
\dfrac{1}{|x_1|^{nq}}\r)dx_1\\
&&\ls {\gz^{n-(n+\ez_1-\az_1)q}}+\gz^{n-nq}.
\end{eqnarray*}
Choosing $\dz\equiv-\max\{n-nq,\,n-(n+\ez_1-\az_1)q\}>0$, we have $G_1\ls\gz^{-\dz}$.

To estimate $G_2$, by the vanishing moments of $a$, we have
\begin{eqnarray*}
&&[b,\,T](a)(x_1,x_2)\\
&&\quad=\dint_{\rnm} [K(x_1,y_1,x_2,y_2)-K(x_1,0,x_2,y_2)
-K(x_1,y_1,x_2,0)+K(x_1,0,x_2,0)]\\
&&\quad\quad\times[b(x_1,x_2)-b(x_1,x_2)
-b(y_1,x_2)+b(y_1,y_2)]a(y_1,y_2)\,dy_1\, dy_2\\
&&\quad\quad +\dint_{\rnm} [K(x_1,y_1,x_2,0)-K(x_1,0,x_2,0)]\\
&&\quad\quad\times[b(x_1,0)-b(x_1,y_2)
-b(y_1,0)+b(y_1,y_2)]a(y_1,y_2)\,dy_1\, dy_2\\
&&\quad\quad +\dint_{\rnm} [K(x_1,0,x_2,y_2)-K(x_1,0,x_2,0)]\\
&&\quad\quad\times[b(0,x_2)-b(y_1,x_2)
-b(0,y_2)+b(y_1,y_2)]a(y_1,y_2)\,dy_1\, dy_2\\
&&\quad\quad + \dint_{\rnm} K(x_1,0,x_2,0)\\
&&\quad\quad\times[b(0,0)-b(y_1,0)
-b(0,y_2)+b(y_1,y_2)]a(y_1,y_2)\,dy_1\, dy_2\\
&&\quad\equiv L_3+L_4+L_5+L_6.
\end{eqnarray*}

Notice that if $x_1\not\in \gz I$ and  $y_1\in I$,  then $|y_1|\le
|x_1|/2$ and $|x_1-y_1|\le 2|x_1|$; if $x_2\not\in \gz_2J$ and
$y_2\in J$,  then $|y_2|\le |x_2|/2$ and $|x_2-y_2|\le 2|x_2|$. Thus
for $x_1\not\in \gz I$ and $x_2\not\in\gz_2 J$, by Definition 1.4,
(K1) through (K4), (R3) and the H\"older inequality, we obtain
$$
|L_3|
\ls\dint_{I}\dint_{J}\dfrac{|y_1|^{\ez_1}}{|x_1|^{n+\ez_1-\az_1}}
\dfrac{| y_2|^{\ez_2}}{|x_2|^{m+\ez_2-\az_2}}|a(y_1,\,y_2)|\,dy_1\, dy_2\\
\ls \dfrac{1}{|x_1|^{n+\ez_1-\az_1}}
\dfrac{1}{|x_2|^{m+\ez_2-\az_2}};
$$

$$
|L_4|\ls\dint_{I}\dint_{J}\dfrac{|y_1|^{\ez_1}}{|x_1|^{n+\ez_1-\az_1}}
\dfrac{| y_2|^{\az_2}}{|x_2|^m}|a(y_1,\,y_2)|\,dy_1\, dy_2
\ls\dfrac{1}{|x_1|^{n+\ez_1-\az_1}} \dfrac{1}{|x_2|^m};
$$
$$
|L_5|\ls\dint_{I}\dint_{J}\dfrac{|y_1|^{\az_1}}{|x_1|^n} \dfrac{|
y_2|^{\ez_2}}{|x_2|^{m+\ez_2-\az_2}}|a(y_1,\,y_2)|\,dy_1\, dy_2 \ls
\dfrac{1}{|x_1|^n} \dfrac{1}{|x_2|^{m+\ez_2-\az_2}};
$$
and
$$
|L_6|\ls\dint_{I}\dint_{J}\dfrac{|y_1|^{\az_1}}{|x_1|^n} \dfrac{|
y_2|^{\az_2}}{|x_2|^m}|a(y_1,\,y_2)|\,dy_1\,
dy_2\ls\dfrac{1}{|x_1|^n} \dfrac{1}{|x_2|^m}.
$$
From this together with $n-(n+\ez_1-\az_1)q<0$, $n-nq<0$,
$m-(m+\ez_2-\az_2)q<0$ and $m-mq<0$, it follows that
\begin{eqnarray*}
G_2\ls&&\dint_{x_1\not\in \gz I}\dint_{x_2\not\in
\gz_2 J}\lf(|L_3|^q+|L_4|^q+|L_5|^q+|L_6|^q\r)\,dx_1\, dx_2\\
&&\ls
\dint_{x_1\not\in \gz I}\dint_{x_2\not\in
\gz_2 J}\lf[\dfrac{1}{|x_1|^{(n+\ez_1-\az_1)q}}
\dfrac{1}{|x_2|^{(m+\ez_2-\az_2)q}}
+\dfrac{1}{|x_1|^{(n+\ez_1-\az_1)q}}
\dfrac{1}{|x_2|^{mq}}\r.\\
&&\quad\lf.+\dfrac{1}{|x_1|^{nq}}
\dfrac{1}{|x_2|^{(m+\ez_2-\az_2)q}}+\dfrac{1}{|x_1|^{nq}}
\dfrac{1}{|x_2|^{mq}}\r]\,dx_1\, dx_2\\
&&\ls\gz^{n-(n+\ez_1-\az_1)q}+\gz^{n-nq}.
\end{eqnarray*}
This shows $G_2\ls\gz^{-\dz}$, which together with
$G_1\ls\gz^{-\dz}$ gives (3.1) and the proof of Theorem 1.3 is
therefore complete.
\end{pfs}

\begin{rem}{\hspace{-0.2cm}}{.}\label{r3.1}\
\rm The restriction $\az_1\le \min\{n/2,\,1\}$ is to guarantee the
boundedness of the commutator $[b,\,T]$ from $L^2(\rnm)$ to
$L^{q_1}(\rnm)$ with $1/q_1=1/p-\az_1/n$; see Theorem 1.2. Since
the $L^2(\rnm)$ norm appears in the definition of $H^p(\rnm)$
rectangular atoms, we need this boundedness of the commutator
$[b,\,T]$ in the proof of Theorem 1.3; see Corollary 1.1.
\end{rem}

\begin{tabular}{ll}
\vspace{-0.3cm}
Der-Chen Chang& Dachun Yang (Corresponding author)\\
\vspace{-0.5cm}
Department of Mathematics&School of Mathematical Sciences\\
\vspace{-0.5cm}
Georgetown University&Beijing Normal University\\
\vspace{-0.5cm}
Washington D. C. 20057&Laboratory of Mathematics and Complex Systems\\
\vspace{-0.5cm}
U. S. A.&Ministry of Education\\
\vspace{-0.5cm}
&Beijing 100875\\
\vspace{-0.2cm}
&People's Republic of China\\
\bigskip
chang@georgetown.edu&dcyang@bnu.edu.cn\\
\end{tabular}

\begin{tabular}{ll}
\vspace{-0.3cm}
Yuan Zhou\\
\vspace{-0.5cm}
School of Mathematical Sciences\\
\vspace{-0.5cm}
Beijing Normal University\\
\vspace{-0.5cm}
Laboratory of Mathematics and Complex Systems\\
\vspace{-0.5cm}
Ministry of Education\\
\vspace{-0.5cm}
Beijing 100875\\
\vspace{-0.2cm}
People's Republic of China\\
 yuanzhou@mail.bnu.edu.cn
\end{tabular}

\section*{References}

\medskip

\begin{enumerate}

\vspace{-0.3cm}
\bibitem[1]{a} T. Aoki, {\it Locally bounded linear topological spaces,}
Proc. Imp. Acad. Tokyo  18 (1942), 588-594.

\vspace{-0.3cm}
\bibitem[2]{b} M. Bownik, {\it Boundedness of operators on
Hardy spaces via atomic decompositions,} Proc. Amer. Math. Soc. 133
(2005), 3535-3542.

\vspace{-0.3cm}
\bibitem[3]{cf1}
S. A. Chang and R. Fefferman, {\it A continuous version of duality
of $H\sp{1}$ with BMO on the bidisc,}
 Ann. of Math. (2) 112 (1980), 179-201.

\vspace{-0.3cm}
\bibitem[4]{cf2}
 S. A. Chang and R. Fefferman,
 {\it The Calder\'on-Zygmund decomposition on product domains,}
Amer. J. Math. 104 (1982), 455-468.

\vspace{-0.3cm}
\bibitem[5]{cf3}
 S. A. Chang and R. Fefferman, {\it Some recent developments in Fourier
analysis and $H\sp p$-theory on product domains,}  Bull. Amer. Math.
Soc. (N. S.) 12 (1985), 1-43.

\vspace{-0.3cm}
\bibitem[6]{chm} W. Chen, Y. Han and C. Miao, {\it Bi-commutators of
fractional integrals on product spaces,} Math. Nachr. 281 (2008),
1108-1118.

\vspace{-0.3cm}
\bibitem[7]{co} R. R. Coifman,
{\it A real variable characterization of $H^p$,} Studia Math. 51
(1974), 269-274.

\vspace{-0.3cm}
\bibitem[8]{f81} R. Fefferman, {\it Singular integrals on
product domains,} Bull. Amer. Math. Soc. (N. S.) 4 (1981), 195-201.

\vspace{-0.3cm}
\bibitem[9]{f85} R. Fefferman, {\it Singular integrals on product $H\sp p$
spaces,} Rev. Mat. Iber. 1:2 (1985), 25-31.

\vspace{-0.3cm}
\bibitem[10]{fe1} R. Fefferman, {\it Calder\'on-Zygmund theory for product domains:
$H\sp p$ spaces,}   Proc. Nat. Acad. Sci. U. S. A. 83 (1986),
840-843.

\vspace{-0.3cm}
\bibitem[11]{fe2}
 R. Fefferman, {\it Harmonic analysis on product spaces,} Ann.
of Math. (2) 126 (1987), 109-130.

\vspace{-0.3cm}
\bibitem[12]{fs} R. Fefferman and E. M. Stein, {\it Singular integrals
on product spaces,} Adv. in Math. 45 (1982), 117-143.

\vspace{-0.3cm}
\bibitem[13]{fl}
 S. H. Ferguson and M. T. Lacey, {\it A characterization
of product BMO by commutators,} Acta Math. 189 (2002), 143-160.

\vspace{-0.3cm}
\bibitem[14]{fjw}
 M. Frazier, B. Jawerth and G. Weiss, {\it Littlewood-Paley Theory and
the Study of Function Spaces,} CBMS Regional Conference Series in
Mathematics, 79, the American Mathematical Society, Providence, R.
I., 1991. \vspace{-0.3cm}

\bibitem[15]{gr}
 J. Garc\'ia-Cuerva and  J. L. Rubio de Francia, {\it Weighted Norm
Inequalities and Related Topics,} North-Holland, 1985.

\vspace{-0.3cm}
\bibitem[16]{g} L. Grafakos, {\it Classical Fourier Analysis},
Second Edition, Graduate Texts in Math., No. 249,
Springer, New York, 2008.

\vspace{-0.3cm}
\bibitem[17]{h89}
Y. Han,  {\it A problem in $H\sp p(\rr\sp 2\sb +\times \rr\sp 2\sb +)$ spaces,}
Chinese Sci. Bull. 34 (1989), 617--622.

\vspace{-0.3cm}
\bibitem[18]{hy}
Y. Han and D. Yang, {\it $H\sp p$ boundedness of Calder\'on-Zygmund
operators on product spaces,} Math. Z. 249 (2005),  869-881.

\vspace{-0.3cm}
\bibitem[19]{j} J.-L. Journ\'e, {\it Calder\'on-Zygmund operators on
product spaces,} Rev. Mat. Iber.  1 (1985), 55-91.

\vspace{-0.3cm}
\bibitem[20]{hyz1} G. Hu, D. Yang and Y. Zhou,
{\it Boundedness  of singular integrals in Hardy spaces on spaces of
homogeneous type,} Taiwanese J. Math. 13 (2009),
91-135.

\vspace{-0.3cm}
\bibitem[21]{La}
 R. H. Latter,  {\it A characterization of $H^p(\rn)$ in terms of
atoms,} Studia Math. 62 (1978), 93-101.

\vspace{-0.3cm}
\bibitem[22]{mtw} Y. Meyer, M. Taibleson and  G. Weiss,
 {\it Some functional analytic properties of the spaces
$B\sb q$ generated by blocks,} Indiana Univ. Math. J. {\bf 34}
(1985), 493-515.

\vspace{-0.3cm}
\bibitem[23]{mc} Y. Meyer and R. R. Coifman,  {\it Wavelets.
Calder\'on-Zygmund and multilinear operators,} Cambridge University
Press, Cambridge, 1997.


\vspace{-0.3cm}
\bibitem[24]{p} J. Pipher,  {\it Journ\'s covering lemma and its extension to
higher dimensions,} Duke Math. J. 53 (1986),  683-690.

\vspace{-0.3cm}
\bibitem[25]{st} E. M. Stein, {\it Singular Integrals and
Differentiability Properties of Functions,} Princeton Univ. Press,
Princeton, N. J., 1970.

\vspace{-0.3cm}
\bibitem[26]{tw}
M. H. Taibleson and G. Weiss, {\it The molecular characterization of
certain Hardy spaces,} Ast\'erisque 77 (1980), 67-149.

\vspace{-0.3cm}
\bibitem[27]{yz} D. Yang and Y. Zhou, {\it Boundedness of
Marcinkiewicz integrals and their commutators in $H^1(\rnm)$,} Sci.
China Ser. A 49 (2006), 770-790.

\vspace{-0.3cm}
\bibitem[28]{yz2} D. Yang and Y. Zhou,
{\it A boundedness criterion via atoms for linear operators in Hardy
spaces,} Constr. Approx. 29 (2009), 207-218.

\vspace{-0.3cm}
\bibitem[29]{yz3} D. Yang and Y. Zhou,
{\it Boundedness of sublinear operators
in Hardy spaces on RD-spaces via atoms,}
J. Math. Anal. Appl. 339 (2008),  622-635.

\vspace{-0.3cm}
\bibitem[30]{y} K. Yabuta, {\it A remark on the $(H^1, L^1)$
boundedness,} Bull. Fac. Sci. Ibaraki Univ. Ser. A 25 (1993), 19-21.
\end{enumerate}
\end{document}